\input ./style/arxiv-general.cfg
\documentclass[aop,MSNbibl,seceqn,dvips]{arximspdf}
\makeatletter
   \@ifpackageloaded{graphicx}{}{\usepackage{graphicx}}
\makeatother
\usepackage{mathrsfs,mathbh}

\doi{10.1214/13-AOP902} 
\volume{43}
\issue{3}
\pubyear{2015}
\firstpage{1419}
\lastpage{1455}

\makeatletter
\newcommand{\rrvert}{\vert}
\newcommand{\llvert}{\vert}
\newcommand{\eqref}[1]{(\ref{#1})}
\newtheorem{thmm}{Theorem}[section]
\newtheorem{lem}[thmm]{Lemma}
\newtheorem{prop}[thmm]{Proposition}

\newproclaim{rem}[thmm]{Remark}
\newproclaim{ex}[thmm]{Example}
\newproclaim{defn}[thmm]{Definition}
\newproclaim{assumptionA}{Assumption}
\newcommand{\varcal}[1]{\mathscr{#1}}
\newcommand{\R}{\mathbb R}
\newcommand{\N}{\mathbb N}
\newcommand{\B}{\varcal{B}}
\newcommand{\x}{\mathbf{x}}
\newcommand{\sgn}{\operatorname{sgn}}
\newcommand{\F}{{\varcal{F}}}
\newcommand{\ind}{\mathbh{1}}
\newcommand{\M}{\mathcal{M}}
\newcommand{\La}{{\Lambda}}
\newcommand{\be}{{\beta}}
\newcommand{\bck}{}
\newcommand{\iy}{{\infty}}
\newcommand{\eps}{{\varepsilon}}
\renewcommand{\P}{\mathbb{P}}
\makeatother

\begin{document}
\begin{frontmatter}

\title{Second-order asymptotics for the block counting process in a
class of regularly varying $\bolds{\Lambda}$-coalescents}
\runtitle{Second-order asymptotics for the block counting process\hspace*{8pt}}

\begin{aug}
\author[A]{\fnms{Vlada} \snm{Limic}\thanksref{T1}\ead[label=e1]{vlada.limic@math.u-psud.fr}}
\and
\author[B]{\fnms{Anna} \snm{Talarczyk}\corref{}\thanksref{T2}\ead[label=e2]{annatal@mimuw.edu.pl}}
\runauthor{V. Limic and A. Talarczyk}
\thankstext{T1}{Supported in part by the ANR MANEGE grant.}
\thankstext{T2}{Supported in part by MNiSzW Grant N N201 397537
and NCN Grant DEC-2012/07/B/ST1/  03417 (Poland).}
\affiliation{Universit\'e Paris-Sud and University of Warsaw}
\address[A]{CNRS UMR 8628\\
Laboratoire de Math\'ematiques\\
Universit\'e Paris-Sud\\
B\^atiment 425\\
91405 Orsay\\
France\\
\printead{e1}} 
\address[B]{Institute of Mathematics\\
University of Warsaw\\
ul.~Banacha 2\\
02-097 Warszawa\\
Poland\\
\printead{e2}}
\end{aug}

\received{\smonth{4} \syear{2013}}
\revised{\smonth{11} \syear{2013}}

%
\begin{abstract}
Consider a standard $\La$-coalescent that comes down from infinity.
Such a coalescent starts from a configuration consisting of infinitely
many blocks at time $0$, but its number of blocks $N_t$
is a finite random variable at each positive time $t$.
Berestycki et al.
[\textit{Ann. Probab.} \textbf{38} (2010) 207--233]
found the first-order approximation $v$ for
the process $N$ at small times.
This is a deterministic function satisfying $N_t/v_t \to1$ as $t\to0$.
The present paper reports on the first progress in the study of the
second-order asymptotics for $N$ at small times.
We show that, if the driving measure $\Lambda$ has a density near zero
which behaves as $x^{-\beta}$ with $\beta\in(0,1)$, then the process
$(\varepsilon^{-1/(1+\beta)}(N_{\varepsilon t}/v_{\varepsilon
t}-1))_{t\ge0}$ converges in law as $\varepsilon\to0$ in the
Skorokhod space to a totally skewed $(1+\beta)$-stable process.
Moreover, this process is a unique solution of a related stochastic
differential equation of Ornstein--Uhlenbeck type, with a completely
asymmetric stable L\'evy noise.

\end{abstract}

%
\begin{keyword}[class=AMS]
\kwd[Primary ]{60J25}
\kwd[; secondary ]{60F17}
\kwd{92D25}
\kwd{60G52}
\kwd{60G55}
\end{keyword}
\begin{keyword}
\kwd{$\La$-coalescent}
\kwd{coming down from infinity}
\kwd{second-order approximations}
\kwd{stable L\'evy process}
\kwd{Ornstein--Uhlenbeck process}
\kwd{Poisson random measure}
\end{keyword}

\pdfkeywords{60J25, 60F17, 92D25, 60G52, 60G55, Lambda-coalescent, coming down from infinity, second-order approximations, stable Levy process, Ornstein-Uhlenbeck process, Poisson random measure}
\end{frontmatter}
%
\section{Introduction and main results}\label{sec1}
\subsection{Background}
\label{S:i}
The $\La$-coalescents were introduced and first studied independently
by Pitman \cite{pit99} and Sagitov \cite{sag99} and were also
considered in a contemporaneous work of Donnelly and Kurtz \cite{dk99}.
They are useful models of genealogical trees of populations that evolve
under the assumption of unbounded variance in the reproduction
(resampling) mechanism.
Berestycki et al. \cite{bbl1} derive the first-order approximation for
the number of blocks in a general standard $\La$-coalescent that comes
down from infinity. The present work initiates the study of the
second-order approximation for the same process.
We next recall the basic definitions, mention some of the landmark
results and present the motivation for the problem we resolved in this work.
For recent overviews of the literature, we refer the reader to \cite
{bertoin,ensaios}.

Let $\La$ be an arbitrary finite measure on $[0,1]$.
We denote by $(\Pi_t, t \ge0)$ the associated $\Lambda$-coalescent.
This Markov jump process $(\Pi_t, t \ge0)$ takes values in the set of
partitions of $\{1, 2, \ldots\}$. Its law is specified by the
requirement that, for any $n\in\N$, the restriction $\Pi^n$ of $\Pi$
to $\{1, \ldots, n\}$ is a continuous-time Markov chain with
the following transitions: whenever $\Pi^n$ has $ b \in[2, n]$ blocks,
any given $k$-tuple of blocks coalesces at rate $\lambda_{b,k}:= \int_{[0,1]} r^{k-2}(1-r)^{b-k} \La(dr)$.
The total mass of $\La$ can be scaled to $1$.
This is convenient for the analysis, and corresponds to a constant time
rescaling of the process. Henceforth, we assume that $\La$ is a
probability measure.

The \emph{standard} $\La$-coalescent starts from the trivial
configuration $\{\{i\}\dvtx i\in\N\}$.
Let us denote by $N^{\La}(t)$ [or $N(t)$ if clear from the context]
the number of blocks of $\Pi(t)$ at time $t$.
If $\P(N^{\La}(t) <\infty, \forall t > 0)=1$, the coalescent is said to
\emph{come down from infinity}.
As part of his thesis work, Schweinsberg \cite{sch1} derived the
following criterion: the (standard) $\La$-coalescent
\emph{comes down from infinity} (CDI) if and only if
%
\begin{equation}
\label{E:cond sch} \sum_{b=2}^\infty \Biggl(
\sum_{k=2}^b (k-1) \pmatrix{b
\cr
k}
\lambda _{b,k} \Biggr)^{-1} <\infty.
\end{equation}
Let
%
\begin{equation}
\label{e:Psi_star} \Psi^*(q)=\int_0^1
\bigl(e^{-yq}-1+qy \bigr)\frac{\Lambda(dy)}{y^2}.
\end{equation}
Bertoin and Le Gall \cite{blg3} obtained an equivalent condition: $\La
$-coalescent CDI if and only if
%
\begin{equation}
\int_a^\infty\frac{1}{\Psi^*(q)}\,dq<\infty\qquad \mbox{for some
(and then all) }a>0. \label{e:cdi}
\end{equation}

Throughout the paper, we will assume \eqref{e:cdi}.
Let $N=(N_t, t\ge0)$ be the block counting process defined above,
so that $N(0)=\infty$ and $\P(N_t<\infty) =1$ for all \mbox{$t>0$}.
As indicated above, in \cite{bbl1}, Theorem~1 it is shown that,
solely under \eqref{e:cdi},
there exists a ``law of large numbers'' approximation for the block
counting process, more precisely,
%
\begin{equation}
\label{Espeed vstar} \lim_{t\to0+} N_t/v_t^*=1\qquad
\mbox{almost surely},
\end{equation}
where $v^*$ is uniquely determined by $\int_{v_t^*}^\iy\frac{1}{\Psi
^*(q)}\,dq = t$, for all $t>0$.
Any function satisfying \eqref{Espeed vstar}
is called a \emph{speed of coming down from infinity}, or a \emph{speed
of CDI}.

Instead of $\Psi^*$ we choose to work with $\Psi\dvtx[1,\infty
)\mapsto
\R
_+$ defined by
%
\begin{equation}
\label{e:Psi} \Psi(q)=\int_0^1
\bigl((1-y)^q-1+qy \bigr)\frac{\Lambda(dy)}{y^2}.
\end{equation}
This function is different from $\Psi$ used in \cite{bbl1} (which is
now our $\Psi^*$).
Moreover, our $\Psi$ appeared as $\bar{\Psi}$ in \cite
{bbl3,me_xi,vl_habi,sampl} where it was already noted that
this function arises from the model in a more natural way [see also
\eqref{e:eqN1} and \eqref{e:Psi2}],
and
it may be more convenient for analysis than $\Psi^*$.
It is not difficult to see that $\Psi$ and $\Psi^*$ have the same
asymptotic behavior at $\infty$ (see Lemma~\ref{lem:3.1} or \cite
{vl_habi,bbl3}), and that therefore~(\ref{E:cond sch}) and (\ref
{e:cdi}) are further equivalent to
%
\begin{equation}
\int_a^\infty \frac{1}{\Psi(q)} \,dq<\infty\qquad
\mbox{ for some (and then all) }a>1. \label{e:cdi2}
\end{equation}
Moreover,
if we define $v\dvtx\R_+\mapsto\R_+$ by
%
\begin{equation}
\label{e:vt} t=\int_{v_t}^\infty \frac{1}{\Psi(q)}
\,dq,
\end{equation}
then (see Lemma~\ref{lem:2.2})
$v_t\sim v_t^*$ as $t\to0$, and so $v$ is also a speed of CDI for the
corresponding $\La$-coalescent.

From the results of Berestycki et al. \cite{bbl1}, it follows that the
asymptotic behavior of the speed $v_t$ of CDI for small $t$ depends
very strongly on the behavior of the driving measure $\Lambda$ near $0$.
This is caused by the fact that the behavior of $\Lambda$ near $0$ is
linked to
the asymptotics of $\Psi(q)$ as $q\to\infty$ by a result of a
tauberian nature. For example, if for small $x$,
%
\begin{equation}
\Lambda(dx)\approx x^{-\beta} \,dx\qquad \mbox{with } \beta\in(0,1),
\label{e:Lambda_informal}
\end{equation}
then $v_t \sim C t^{-{1}/\beta}$, for some $C\in(0,\iy)$, as
$t\to0$ (see Lemma~\ref{lem:3.2}). Note that \eqref{e:Lambda_informal}
is understood in the sense of assumption (\ref{e:A}) in Section~\ref{S:result}.

A natural question is to study the second-order fluctuations of $N$
about its speed of CDI. In particular, one wishes to
understand how close is
$\frac{N_t}{v_t}$ to $1$ at small times, and if this proximity can be
measured in some regular (and universal) way.
In the present paper, we address this problem by considering the
fluctuations in a functional sense, with time scaled by $\varepsilon
\to0$.
More precisely, we investigate the convergence in law of the processes
%
\begin{equation}
\label{e:fluct} \biggl( r(\eps) \biggl( \frac{N_{ \eps t}}{v_{ \eps t}} - 1 \biggr), t \geq0
\biggr),
\end{equation}
were $r(\varepsilon)$ is an appropriately chosen normalization so that
the limit process is nontrivial.

It turns out that both the normalization $r(\varepsilon)$ and the limit
process again depend on the behavior of $\Lambda$ near $0$. The
singularity exponent $\be$ of the density\vspace*{1pt} of $\La$ near $0$ decides the
rate of convergence
of $\frac{N_t}{v_t}$ and, therefore, of $\frac{N_t}{v_t^*}$, to 1.

\subsection{Main results}
\label{S:result}
We assume that the coalescent does not have a Kingman part and also that
$\Lambda(\{1\})=0$, so that the $\La$-coalescent either comes down from
infinity or stays infinite forever
(see Pitman \cite{pit99}).
We formalize \eqref{e:Lambda_informal} in the following way, making it
our main assumption.

\begin{assumptionA*}
$\Lambda(\{0\})=\Lambda(\{1\})=0$. Moreover,
there exists $y_0\le1$ such that
%
{\renewcommand{\theequation}{A}
\begin{equation}
\label{e:A} 
\La(dy)=g(y)\,dy,\qquad y \in[0,y_0]\quad \mbox{and}
\quad\lim_{y\to0+}g(y)y^\beta=A
\end{equation}}
\hspace*{-3pt}for some $0<\beta<1$ and $0<A<\infty$.
\end{assumptionA*}

\begin{rem}
\label{rem:2.3}
(a) Condition $\beta>0$ ensures that the $\Lambda$ coalescent satisfies~\eqref{e:cdi2}, hence that it comes down from infinity, since it is not
difficult to see that \eqref{e:A} implies that $\Psi(q)\sim C
q^{1+\beta
}$ as $q\to\infty$ (see also Lemma~\ref{lem:3.2} below).
Condition $\beta<1$ is clear, since $\Lambda$ has to be a finite
measure.

(b) Assumption \eqref{e:A} is satisfied by all the Beta-coalescents
that come down from infinity, that is, all the coalescents where
$\Lambda$ has density of the form
$g(y)=\frac{1}{B(1-\beta,a)} y^{-\beta}(1-y)^{a-1}$, for some $0<\beta
<1$ and $a>0$ and the normalizing constant is the appropriately
evaluated Beta function.

(c) By Lemma~\ref{lem:3.1} in the next section, $\Psi$ is a continuous
and strictly increasing function on $[1,\iy)$, strictly positive on
$(1,\iy)$, and $\int_1^\iy \,dq/\Psi(q)\geq\int_1^\iy \,dq/q(q-1)=\iy$.
This, together with CDI, implies that $v$ given by \eqref{e:vt}
is a well defined strictly decreasing function on $(0,\infty)$ and it
takes values in $(1,\infty)$.

Further properties of $v$ and $\Psi$ can be found in Section~\ref
{S:preliminary}.
Under assumption \eqref{e:A}, we can obtain precise
asymptotics of the speed of coming down from infinity $v$ and the
function $\Psi$; see Lemma~\ref{lem:3.2}. In particular, as $t\to0$ we
have $ v_t\sim v^*_t\sim K_1 t^{-{1}/\beta}$,
where
%
\setcounter{equation}{9}
\begin{equation}
K_1= \biggl(\frac{1+\beta}{A \Gamma(1-\beta)} \biggr)^{{1}/\beta}, \label{e:K1}
\end{equation}
and where $\Gamma$ is the Gamma function.
\end{rem}

We shall study the asymptotic behavior, as $\varepsilon\to0$, of the
process $X_\varepsilon=(X_\varepsilon(t))_{t\ge0}$ defined by
%
\begin{equation}
\label{e:Xe} X_\varepsilon(0)=0 \quad\mbox{and}\quad X_\varepsilon(t)=
\varepsilon ^{-{1}/{(1+\beta)
}} \biggl(\frac{N_{\varepsilon t}}{v_{\varepsilon t}} -1 \biggr), \qquad t>0.
\end{equation}
For each $B\in\B(\R)$ Borel set, let $\vert B\vert$ denote its
Lebesgue measure.
Let $\M$ be an independently scattered $(1+\beta)$-stable random
measure on $\R$ with skewness intensity $1$.
That is, for each $B\in\B(\R)$ such that $0<\vert B\vert<\infty$,
$\M(B)$
is a $(1+\beta)$-stable random variable
with characteristic function
\[
\exp \biggl\{-\vert B\vert\vert z\vert^{1+\beta} \biggl(1-i(\sgn z)\tan
\frac
{\pi
(1+\beta)}2 \biggr) \biggr\}, \qquad z\in\R,
\]
$\M(B_1), \M(B_2),\ldots$ are independent whenever $B_1, B_2, \ldots$
are disjoint sets,
and $\M$ is $\sigma$-additive a.s. (see Samorodnitsky and Taqqu \cite
{ST}, Definition~3.3.1).

We are now ready to state the main result.

\begin{thmm}
\label{thmm:main}
Assuming \eqref{e:A}, the process $X_\varepsilon$ defined in \eqref{e:Xe}
converges in law in the Skorokhod space $D([0,\infty))$ equipped with
$J_1$ topology to a $(1+\beta)$-stable process $Z=(Z_t)_{t\ge0}$
given by
%
\begin{equation}
\label{e:Z} Z(t)=-\frac{K} t\int_0^tu
\M(du),\qquad t>0, Z(0)=0,
\end{equation}
where $K$ is the following positive constant:
%
\begin{equation}
\label{e:K} K= \biggl(-A \int_0^\infty
\bigl(e^{-y}-1+y \bigr)y^{-2-\beta}\,dy \cos\frac{\pi
(1+\beta)}2
\biggr)^{{1}/{(1+\beta)}}.
\end{equation}
\end{thmm}

The proof of this theorem is given in Section~\ref{S:pf_main}.

\begin{rem}
\label{rem:2.5}
(a) The integral in \eqref{e:Z} is understood in the sense of Chapter~3 of \cite{ST}.

(b) The process $Z$ can be also expressed as
\[
Z(t)=-\frac{K}t\int_0^t u
\,dL_u,\qquad t>0, Z(0)=0,
\]
where $L$ is the $(1+\beta)$-stable
totally skewed to the right (having no negative jumps) L\'evy process.
Moreover, $Z$ solves the following stochastic differential equation of
the Ornstein--Uhlenbeck type:
%
\begin{equation}
\label{e:eqZ} Z(t)=-\int_0^t
s^{-1} Z(s) \,ds-KL(t).
\end{equation}

(c) It was already mentioned (cf. Remark~\ref{rem:2.3}) that assumption
(\ref{e:A}) is satisfied by Beta-coalescents which come down from infinity.
Theorem~\ref{thmm:main} shows that, from the point of view of behavior
of $N_t$, $v_t$ and $N_t/v_t-1$ near $0$, any $\Lambda$-coalescent
satisfying (\ref{e:A}) resembles a corresponding Beta-coalescent (or rather a
class of Beta-coalescents) having driving measure(s) of the form
$\operatorname{Beta}(1-\beta, a)$, for some $a>0$.

The fact that the limit process is $(1+\beta)$-stable can be explained
by observing that for each $\beta\in(0,1)$, one member of the above
family [notably the $\operatorname{Beta}(1-\beta,1+\beta)$-coalescent] was obtained
from genealogies of populations with supercritical infinite variance
branching both by Sagitov \cite{sag99} [in his setting, the branching
mechanism has generating function
$1-\frac{1+\beta}{\beta}(1-s)+\frac{1}\beta(1-s)^{1+\beta}$] and by
Schweinsberg \cite{Schweinsberg2003} (in his setting, the probability
that the individual has $k$ or more offspring decays like $k^{-(1+\beta)}$).
It is well known that branching laws of this type are in the domain of
attraction of the $(1+\beta)$-stable law. Moreover, the limits of
fluctuations related to infinite variance branching systems of type
$1+\beta$ are usually $(1+\beta)$-stable. (See, e.g., Iscoe
\cite
{Iscoe} Theorem~5.4 and 5.6 and
Bojdecki et al. \cite{functlim3}.)
Another connection is due to \cite{7}, relating $\operatorname{Beta}(1-\beta,1+\beta
)$-coalescents to continuous state $(1+\beta)$-stable processes.
The limit process is naturally totally skewed to the left, as $N_t$
only has negative jumps, hence so does $X_{\varepsilon}$.
\end{rem}

We also wish to mention here a related work of Schweinsberg \cite
{Schweinsberg2012}, where fluctuations of the number of blocks of the
Bolthausen--Sznitman coalescent were investigated (see Theorem~1.7 in
\cite{Schweinsberg2012}). This is a different setting from ours, since
the Bolthausen--Sznitman coalescent does not come down from infinity
[\eqref{e:Lambda_informal} holds in this case with $\beta=0$].
Schweinsberg investigated appropriately rescaled fluctuations of the
number of blocks of the Bolthausen--Sznitman coalescent starting from
$n$ blocks in the limit as $n\to\infty$. It is interesting to note
that the limit in \cite{Schweinsberg2012} involves a totally skewed
$1$-stable process.

Another interesting fact is that the present analysis (in the sense of
functional convergence) has not been carried out even for the
case of the Kingman coalescent, where $\La$ is the Dirac measure at 0.
It is known in this case that the law of $t^{-{1}/2}(N_t/v_t-1)$
converges to a Gaussian law; see, for example, Aldous \cite{aldous_survey}.
Here, we assume that $\Lambda(\{0\})=0$, so that the coalescent does
not have the Kingman part.
We postpone the study of the complementary setting to a future work.
We conjecture that in the case of the pure Kingman coalescent
(i.e., $\La$ is the Dirac mass at $0$) the limit process in \eqref
{e:fluct} will have a form similar to \eqref{e:Z}, where the
integration with respect to the stable random
measure is replaced by integration with respect to Brownian motion. The
Kingman case, although seemingly easier, cannot be done with our
present technique, since here we rely heavily on the Poisson process
construction of $\Lambda$ coalescents, which is particularly nice if
$\Lambda(\{0\})=0$.

Under assumption \eqref{e:A}, we have $v_t\sim v_t^*\sim w_t=K_1
t^{-{1}/{\beta}}$ (see Lemma~\ref{lem:3.2}). It is therefore
natural to ask whether one obtains the same results if in \eqref{e:Xe}
$v$ is replaced by $v^*$ or $w$.
The answer is positive for $v^*$. For $w^*$, one has to assume
additional regularity of $\Lambda$ near $0$.

Define
$X_\varepsilon^*(0)=0$, $X_\varepsilon^\be(0)=0$ and
%
\begin{eqnarray}
\label{e:Xe star} X_\varepsilon^*(t)&=&\varepsilon^{-{1}/{(1+\beta)}} \biggl(
\frac
{N_{\varepsilon
t}}{v_{\varepsilon t}^*} -1 \biggr),
\nonumber
\\[-8pt]
\\[-8pt]
\nonumber
 X_\varepsilon^\be(t)&=&\varepsilon
^{-{1}/{(1+\beta)}} \biggl((\eps t)^{1/\be}
\frac{N_{\varepsilon t}} {K_1} -1 \biggr),\qquad t>0,
\end{eqnarray}
where $K_1$ is the constant given by \eqref{e:K1}. Let
$\Longrightarrow
$ denote the convergence in law of processes with respect to the
Skorokhod topology.

As a corollary to Theorem~\ref{thmm:main}, we obtain the following results.

%
\begin{thmm}
\label{thmm:main a}
Assume \eqref{e:A}, and let $Z$ and $K$ be as in Theorem~\ref
{thmm:main}. Then
\begin{longlist}[(a)]
\item[(a)] $X_\varepsilon^* \Longrightarrow Z$,
\item[(b)]if moreover $(y^\be g(y)-A) = O(y^\alpha)$, as $y\to0$, for some
$\alpha>\be/(1+\be)$,
then
\[
X_\varepsilon^\be\Longrightarrow Z.
\]
\end{longlist}
\end{thmm}

The
proof is postponed until Section~\ref{S:robust}.

\begin{rem}
\label{rem:natural}
As a counterpart to part (b) in Section~\ref{S:robust b}, we exhibit a
family of counterexamples,
for which $y \mapsto y^\be g(y)$ is not sufficiently H\"older
continuous at $0$, and the
above ``natural extension'' of convergence in Theorem~\ref{thmm:main a}(b)
fails.
In turns out that one does not have to search hard for counterexamples:
the first guess $g(y)= y^{-\be} + y^{\alpha-\be}$, where $\alpha$ is
such that $\alpha<\be/(\be+1)$, already does the trick.
This illustrates a remarkable sensitivity of the second-order
approximation for $N$ with respect to the smoothness of $\La$ near $0$.
\end{rem}

\subsection{Main tools}
\label{ss:main_tools}
When $\Lambda(\{0\})=0$,
one can construct a realization of the corresponding
$\Lambda$-coalescent from a Poisson point process in
the following (now standard) way. Let
%
\begin{equation}
\label{DPPPpi} \pi(\cdot) = \sum_{i \in\N}
\delta_{(T_i,Y_i)}(\cdot)
\end{equation}
be a Poisson point process on $\R_+ \times(0,1)$ with intensity measure
$dt \otimes\nu(dy)$ where $\nu(dy)= y^{-2}\Lambda(dy)$.
Each atom $(t,y)$ of $\pi$ impacts the evolution of $\Pi$ as
follows:
for each block of $\Pi(t-)$ a coin is flipped with probability of
heads equal to $y$; all the blocks corresponding to
coins that come up ``head'' are merged immediately
into one single block, and all the other blocks remain unchanged.
In order to make this construction rigorous,
one initially considers the restrictions
$(\Pi^{(n)}(t),t\ge0)$, since the measure $\nu$ may be infinite (see,
e.g., \cite{ensaios,bertoin}).

Our technique is based on a novel approach, using an explicit
representation of the block counting process in terms of an enriched
Poisson random measure $\pi^E$. This measure $\pi^E$ is defined on a
larger space in such a way that it also includes the information on
(individual block) coloring. One can then write an integral equation
for the number of blocks $N_t$ involving an integral with respect to
$\pi^E$. This equation turns out to be analytically tractable. In our
approach, we rely on the properties of integrals with respect to
Poisson, compensated Poisson and stable random measures, Laplace
transforms of Poisson integrals and of totally skewed stable random
variables, as well as standard tools in the analysis of processes in
the Skorokhod space, for example, the Aldous criterion for tightness.
Moreover, a deterministic lemma from \cite{bbl1}, for comparing
solutions to two different Cauchy (or Cauchy-like) problems, turns out
to be very useful.

The remainder of the paper is organized as follows.
In Section~\ref{S:preliminary}, we give some basic information on the
properties
of $\Psi$ and $v$; in Section~\ref{S:Basic}, we develop the integral
equations for $N$ and $N/v$
and study their basic properties. This is done in a fairly general setting;
in Section~\ref{S:pf_main}, we give the proof of the main result---Theorem~\ref{thmm:main}; in Section~\ref{S:robust}, we prove Theorem~\ref{thmm:main a} and discuss the problem of robustness.

Throughout the paper, $C, C_1, C_2, \ldots$ always denote positive
constants which may be different from line to line.

\section{Preliminary results}
\label{S:preliminary}
In this section, we collect some of the basic properties of $\Psi$ and
$v$ and their relation to the block counting process $N$.
Unless otherwise stated, the facts presented in this section do not
require \eqref{e:A} and are derived for general $\Lambda$.

Recall that $\Psi$ and $v$ are defined by \eqref{e:Psi} and \eqref
{e:vt}, respectively.
Let us also define
%
\begin{equation}
h(q):=\frac{\Psi(q)}{q}. \label{e:h}
\end{equation}
For $0<a\le1$, let $\Psi_a$ (resp., $\Psi^*_a$) be defined by \eqref{e:Psi}
[resp., \eqref{e:Psi_star}]
with $\Lambda(dy)$ replaced by $\Lambda_a(dy)=\ind_{[0,a]}(y)\Lambda(dy)$.

The first lemma concerns the most general setting, up to time-change.

%
\begin{lem}
\label{lem:3.1} Let $\Lambda$ be an arbitrary probability measure on
$[0,1]$ satisfying
$\Lambda(\{0\})=\Lambda(\{1\})=0$.
Then the function $\Psi$ given by \eqref{e:Psi} is well defined on
$[1,\infty)$. In addition,
\begin{longlist}[(iii)]
\item[(i)] $\Psi$ is continuous on $[1,\infty)$ and strictly positive on
$(1,\infty)$,

\item[(ii)] for any $q\ge1$
%
\begin{eqnarray}
\label{e:3.2} \Psi(q)&\le& q(q-1),
\\
\label{e:3.3} 0&\leq&\Psi^*(q)- \Psi(q)\le\frac{q}2,
\end{eqnarray}
\item[(iii)] for any $q\ge1$ and $a \in(0,1)$
%
\begin{eqnarray}
0&\le&\Psi(q)-\Psi_a(q)\le \frac{q}{a},\label{e:3.3a}
\\
0&\le&\Psi^*(q)-\Psi^*_a(q)\le \frac{q}{a}, \label{e:3.8}
\end{eqnarray}
\item[(iv)]
and both $\Psi$ and $h$ are strictly increasing on $[1,\infty)$ and
differentiable on $(1,\infty)$.
\end{longlist}
\end{lem}

Most of these facts are known in the literature but for the benefit of
the reader
we will include a short proof.
Note that \eqref{e:3.3} implies the equivalence of \eqref{e:cdi} and~\eqref{e:cdi2}.

\begin{pf*}{Proof of Lemma \ref{lem:3.1}}
We start with some useful representations for $\Psi$. Clearly, $\Psi
(1)=0$ and if $q>1$ we have
%
\begin{eqnarray}
\label{e:3.4} \Psi(q)&=&q\int_0^1\int
_0^y \bigl(1-(1-r)^{q-1} \bigr)\,dr
\frac
{\Lambda
(dy)}{y^2}
\\
\label{e:3.5}& =& q(q-1)\int_0^1\int
_0^y\int_0^r(1-u)^{q-2}
\,du \,dr\frac{\Lambda
(dy)}{y^2}
\\
\label{e:3.6} &=& q(q-1)\int_0^1\int
_0^1\int_0^r(1-uy)^{q-2}
\,du \,dr{\Lambda(dy)}.
\end{eqnarray}
Representation \eqref{e:3.6} shows that $\Psi$ is finite, continuous on
$[1,\infty)$, and
strictly positive on $(1,\infty)$.
Note that if $q\geq2$, then the integrand in \eqref{e:3.6} is smaller
than $1$ so $\Psi(q) \leq q(q-1)/2$.
The general estimate \eqref{e:3.2}
follows from \eqref{e:3.6}, the fact that for $0\le u,y\le1$ and
$q\ge
1$ we have $(1-uy)^{q-2}\le(1-u)^{-1}$ (easy for $q=1$, and then use
monotonicity) and finally the identity
$\int_0^1 \log(1-r) \,dr=-1$.
The estimates of type \eqref{e:3.3} were already derived in \cite
{bbl1,sampl,vl_habi}.
The lower bound is a consequence of \eqref{e:Psi_star}, \eqref{e:Psi}
and the trivial inequality $(1-y)^q\le e^{-qy}$ for $0\le y\le1$.
The upper bound can obtained, for example, by using \eqref{e:3.4} and
its analogue for $\Psi^*$ that yield
\[
\Psi^*(q)-\Psi(q)=q\int_0^1\int
_0^y \bigl((1-r)^{q-1}-e^{-qr}
\bigr)\,dr\frac{\Lambda(dy)}{y^2},
\]
and observing that $(1-r)^{q-1}-e^{-qr}\le(1-r)^{q-1}-(1-r)^q\le r$
for $0\le r\le1$ and $q\ge1$.
The bound
\eqref{e:3.3a} follows easily from \eqref{e:3.4}, and \eqref{e:3.8} can
be proved via a similar representation for $\Psi^*$.
For (iv), it clearly suffices to show that $h$ is increasing and
differentiable.
This can be easily seen from \eqref{e:3.4}.
\end{pf*}

From now on, we assume that $\Lambda(\{0\})=\Lambda(\{1\})=0$ and that
the $\Lambda$-coalescent comes down from infinity, which is equivalent
to any of \eqref{E:cond sch}, \eqref{e:cdi}, \eqref{e:cdi2}.
By Lemma~\ref{lem:3.1}, $\Psi$ is a continuous and strictly increasing
function on $[1,\iy)$, strictly positive on
$(1,\iy)$ and $\int_1^\iy \,dq/\Psi(q)\geq\int_1^\iy \,dq/q(q-1)=\iy$.
As already mentioned in the \hyperref[sec1]{Introduction}, this
implies that $v$ is a
well defined strictly decreasing function on $(0,\infty)$.
Moreover, $v$ has the following properties.

%
\begin{lem}
\label{lem:2.2}
\textup{(i)} $v_t>1$ for all $t> 0$, $\lim_{t\to0+} v_t = \iy$ and $\lim_{t\to
\iy} v_t =1$,
\begin{longlist}[(iii)]
\item[(ii)] $v$ is differentiable and
%
\begin{equation}
v_t'=-\Psi(v_t), \label{e:2.3a}
\end{equation}
\item[(iii)] in addition
%
\begin{equation}
\lim_{t\to0+}\frac{v_t}{v_t^*}=1. \label{e:2.4}
\end{equation}
\item[(iv)] Therefore,
%
\begin{equation}
\lim_{t\to0+}\frac{N_t}{v_t}=1\qquad \mbox{almost surely,} \label{e:2.5}
\end{equation}
\item[(v)] and for any $p>0$,
%
\begin{equation}
\label{e:2.6} \lim_{t\to0+} E\sup_{0<s\le t}\biggl\vert
\frac{N_s}{v_s}-1\biggr\vert^p=0.
\end{equation}
Moreover, for any $p>0$ there exists $C(p)>0$ such that
%
\begin{equation}
\label{e:2.8} E \sup_{s \ge0} \biggl(\frac{N_s}{v_s}
\biggr)^p\le C(p).
\end{equation}
\end{longlist}
\end{lem}

\begin{rem}
\label{R:pmom vvstar}
Parts (iv) and (v) of Lemma~\ref{lem:2.2} say that $\frac{N_t}{v_t}$
converges to $1$ almost surely and in $L^p$, for any $p>0$.
This was shown with $v^*$ in place of $v$ in \cite{bbl1} Theorems 1
and 2.
Moreover,
in the same article \eqref{e:2.8} was derived, again with $v^*$ in
place of $v$.
(Note that \cite{bbl1} Theorem~2 assumes that $p\geq1$, but this can
be easily extended to all $p\in(0,1)$ by Jensen's inequality.)
Due to \eqref{e:2.4}, one obtains (iv)--(v) without any additional work.
In comparison, Lemma~\ref{lem:p_moment} stated at the end of Section~\ref{S:Basic} is a novel and stronger estimate, important for our analysis.
\end{rem}

We recall next the following elementary estimate that
will be used frequently in the proofs (see \cite{bbl1}, Lemma~10 for
derivation).

%
\begin{lem}
\label{L:bbl1}
Suppose $f,g\dvtx[a,b]\mapsto{\mathbb R}$ are c\`adl\`ag functions
such that
%
\begin{equation}
\label{Ecoupledfun} \sup_{x\in[a,b]} \biggl\llvert f(x) + \int
_a^x g(u) \,du \biggr\rrvert \leq c
\end{equation}
for some $c<\infty$.
If in addition
$f(x)g(x) > 0$, $x\in[a,b]$ whenever $f(x) \neq0$,
then
\[
\sup_{x\in[a,b]} \biggl\llvert \int_a^x
g(u) \,du \biggr\rrvert \leq c \quad\mbox{and}\quad \sup_{x\in[a,b]}\bigl |f(x)\bigr|
\leq2c.
\]
\end{lem}

\begin{pf*}{Proof of Lemma~\ref{lem:2.2}}
We have $\Psi(1)=0$. Moreover, \eqref{e:3.2} shows that $\int_1^\iy
\,dq/\Psi(q)=\infty$. Together with
the strict positivity of $\Psi$ on $(1,\iy)$ and (\ref{e:cdi2}), this
implies that $x\to F(x):=
\int_x^\iy \,dq/\Psi(q)$ maps $(1,\infty)$ bijectively to $(0,\infty)$.
Since $v$ is the inverse of $F$, it is clearly a strictly decreasing
function and (i) holds.
Property (ii) is clear by the definition of $v$ and fundamental
theorem of calculus.
Provided we show the claim in (iii), (iv) is clearly true due to \eqref
{Espeed vstar}.
Similarly,
\[
\frac{N_t}{v_t} -1= \frac{v_t^*}{v_t} \biggl(\frac{N_t}{v_t^*} -1 \biggr) +
\frac{v_t^*}{v_t} -1,
\]
so (iii) and \cite{bbl1} Theorem~2 together imply \eqref{e:2.6}.
The estimate in \eqref{e:2.8} follows easily from \eqref{e:2.6} by the
triangle inequality,
the (decreasing) monotonicity of $N$, and the fact that $v_t\in(1,\iy
)$ for each $t>0$.

In the rest of the argument, we prove (iii). This deterministic
argument is a simplified version of the stochastic (martingale based)
argument for \cite{bbl1}, Theorem~1.
We will show a somewhat stronger statement: $\log\frac
{v_t}{v_t^*}=O(t)$ as $t\to0+$.
In order to do this, for $n\in\N$, $n>1$ define the functions
$v^{(n)}$ and $v^{*,(n)}$ by
\[
t=\int_{v^{(n)}_t}^n \frac{1}{\Psi(q)}\,dq \quad\mbox{and}\quad
t=\int_{v^{*,(n)}_t}^n \frac{1}{\Psi^*(q)}\,dq.
\]
By Lemma~\ref{lem:3.1}, $\Psi$ is strictly positive on $(1,\infty)$ and
it satisfies $\int_1^n\frac{dq}{\Psi(q)}=\infty$, hence $v_t^{(n)}$ is
well defined.
Similarly, it is easy to see (and checked in \cite{bbl1}) that $\Psi^*$
is strictly positive on $(0,\infty)$ and $\int_0^n\frac{dq}{\Psi
^*(q)}=\infty$, so $v^{*,(n)}_t$ is also well defined. Moreover, by
\eqref{e:cdi} and \eqref{e:cdi2} for each $t>0$, we have that
$v^{(n)}_t\nearrow v_t$ and $v^{*,(n)}_t\nearrow v^*_t$ as $n\to\infty
$. The functions $v^{(n)}$ and $v^{*,(n)}$ satisfy equations
\[
v_t^{(n)}=n-\int_0^t
\Psi \bigl(v_s^{(n)} \bigr)\,ds \quad\mbox{and}\quad
v_t^{*,(n)}=n-\int_0^t
\Psi^* \bigl(v_s^{*,(n)} \bigr)\,ds.
\]
Hence, $d\log{v_t^{(n)}}=-\Psi(v_t^{(n)})/v_t^{(n)} \,dt$ and $d\log
{v_t^{*,(n)}}=-\Psi^*(v_t^{*,(n)})/v_t^{*,(n)} \,dt$. This implies that
\[
\log\frac{v_t^{(n)}}{v_t^{*,(n)}} + \int_0^t \biggl[
\frac{\Psi
(v_s^{(n)})}{v_s^{(n)}} - \frac{\Psi
^*(v_s^{*,(n)})}{v_s^{*,(n)}} \biggr] \,ds=0.
\]
Observe also that if $t$ is sufficiently small, then $v^*_t\ge2$.
Hence, there exists a $t_2^*>0$ such that for all sufficiently large
$n$ we have $\inf_{t\in[0,t_2^*]}v_t^{*,(n)}>1$.
For such $n$ and $t\leq t_2^*$, one can rewrite the last identity as
%
\begin{eqnarray}
\label{E:last id} &&\log\frac{v_t^{(n)}}{v_t^{*,(n)}}+ \int_0^t
\biggl[\frac{\Psi
(v_s^{(n)})}{v_s^{(n)}} - \frac{\Psi(v_s^{*,(n)})}{v_s^{*,(n)}} \biggr] \,ds
\nonumber
\\[-8pt]
\\[-8pt]
\nonumber
&&\qquad=\int
_0^t \frac{\Psi^*(v_s^{*,(n)})- \Psi
(v_s^{*,(n)})}{v_s^{*,(n)}} \,ds.
\end{eqnarray}
By \eqref{e:3.3}, the absolute value of the integral on the right-hand
side of this equation is bounded by $\frac{t}2$. Moreover,
by Lemma~\ref{lem:3.1}(iv), the function $q\mapsto\Psi(q)/q$ is
strictly increasing, so we can apply Lemma~\ref{L:bbl1} obtaining
$ \vert\log(v_t^{(n)}/\break v_t^{*,(n)})\vert\le t$. Letting $n\to\infty
$, we get
\begin{equation}
\biggl\vert\log\frac{v_t}{v_t^{*}}\biggr\vert\le t, \label{e:5.4a}
\end{equation}
thus completing the proof.
\end{pf*} 

Under assumption \eqref{e:A}, it is possible to study the asymptotics
of $\Psi$ and $v$ in much more detail, as given by the following lemma.

%
\begin{lem}
\label{lem:3.2}
Assume \eqref{e:A}.
Then
\begin{longlist}[(iii)]
\item[(i)]
%
\begin{equation}
\lim_{q\to\infty}\frac{\Psi(q)}{q^{1+\beta}}=\lim_{q\to\infty
}
\frac
{\Psi^*(q)}{q^{1+\beta}} =\frac{A\Gamma(1-\beta)}{\beta(\beta+1)}, \label{e:3.9}
\end{equation}
\item[(ii)]
%
\begin{equation}
\lim_{t\to0+} tv_t^\beta=\lim
_{t\to0+} t \bigl(v^*_t \bigr)^\beta=
\frac
{1+\beta
}{A \Gamma(1-\beta)}. \label{e:3.10}
\end{equation}
Moreover, there exist $C_1,C_2>0$ such that for all $t>0$
%
\begin{equation}
\label{e:3.10a} C_1 \bigl(t^{-{1}/\beta}\vee1 \bigr)\le
v_t \le C_2 \bigl(t^{-{1}/\beta}\vee1 \bigr).
\end{equation}
\item[(iii)] For $h$ defined by \eqref{e:h}, we have
%
\begin{equation}
\lim_{q\to\infty} q^{1-\beta} h'(q)=
\frac{A \Gamma(1-\beta
)}{1+\beta}, \label{e:3.11}
\end{equation}
moreover,
%
\begin{equation}
\label{e:3.12} \sup_{q\ge1}q^{1-\beta}h'(q)<
\infty.
\end{equation}
\end{longlist}
\end{lem}

\begin{pf}
(i) From assumption \eqref{e:A}, it follows that there exists
$0<a<\frac{1}2$ such that $\Lambda$ has a density $g$ on $[0,a]$ and
%
\begin{equation}
\label{e:A2} \frac{A}2\le\inf_{0<y\le a}
g(y)y^\beta\le\sup_{0<y\le a} g(y)y^\beta
\le2A.
\end{equation}
Due to \eqref{e:3.3}--\eqref{e:3.8}, it suffices to prove \eqref{e:3.9}
with $\Psi^*_a$.
It is immediate to check that $\Psi^*_a(q)=q^2\int_0^1\int_0^1\int_0^r
e^{-quy}\,du\,dr\Lambda_a(dy)$ [note that this is an analogue of \eqref
{e:3.6}]. Hence,
\begin{eqnarray*}
\lim_{q\to\infty}\frac{\Psi^*_a(q)}{q^{1+\beta}}& =&\lim_{q\to\infty}{q^{1-\beta}}
\int_0^1 \int_0^r
\int_0^a e^{-qyu}g(y)\,dy\,du\,dr
\\
&=&\lim_{q\to\infty} \int_0^1\int
_0^r \int_0^{auq}u^{\beta-1}
e^{-y} y^{-\beta} {g \biggl(\frac{y}{qu} \biggr) \biggl(
\frac{y}{qu} \biggr)^\beta} \,dy \,du \,dr
\\
&=&\frac{A\Gamma(1-\beta)}{\beta(1+\beta)},
\end{eqnarray*}
where the second equality is obtained via the substitution $y'=uqy$
(then $y'$ is renamed $y$)
while the third follows by \eqref{e:A}, \eqref{e:A2} and the dominated
convergence theorem.

(ii) Due to \eqref{e:vt} and the fact that $v$ diverges to $\iy$ at
$0$, we have
\[
\lim_{t\to0} tv_t^\beta=\lim
_{x\to\infty} x^\beta\int_x^\infty
\frac
{1}{\Psi(q)}\,dq,
\]
and by the l'Hospital rule and \eqref{e:3.9} we obtain that $ \lim_{t\to0} tv_t^\beta=\frac{1+\beta}{A \Gamma(1-\beta)}$. The same is
true for $v^*$. Finally, note that \eqref{e:3.10a} follows from \eqref
{e:3.10},
the (decreasing) monotonicity of $v$ and the fact that $v_t>1$ for all $t$.

(iii) Let $a$ be as in the proof of part (i). By
\eqref{e:3.4}, we have that
%
\begin{equation}
h=h_a+\tilde h_a, \label{e:3.15}
\end{equation}
where
%
\begin{eqnarray}
h_a(q)&=&\int_0^a \int
_0^y \bigl(1- (1-r)^{q-1} \bigr)\,dr
\frac
{\Lambda(dy)}{y^2}, \label{e:3.15a}
\\
\tilde h_a(q)&=&\int_a^1 \int
_0^y \bigl(1- (1-r)^{q-1} \bigr)\,dr
\frac
{\Lambda(dy)}{y^2}. \label{e:3.15b}
\end{eqnarray}
Then
%
\begin{equation}
\label{e:3.16} h_a'(q)=\int_0^a
\int_0^y \bigl(-\ln(1-r) \bigr)
(1-r)^{q-1} \,dr \frac
{g(y)}{y^2}\,dy
\end{equation}
and
%
\begin{equation}
\label{e:3.17} \tilde h_a'(q)=\int
_a^1\int_0^y
\bigl(-\ln(1-r) \bigr) (1-r)^{q-1}\,dr \frac{\Lambda(dy)}{y^2}.
\end{equation}
In the above expression for $\tilde h_a'$, we substitute $r'=-\ln(1-r)$
and use the obvious estimates to get
%
\begin{equation}
\label{e:3.18} \tilde h_a'(q)\le\frac{1}{a^2}\int
_0^\infty re^{-rq} \,dr=
\frac{1}{a^2q^2}.
\end{equation}
For $h_a'$, we first use the substitution $r'=\frac{r}{y}$ and then
$y'=y(q-1)r'$ to obtain
%
\begin{eqnarray}
\label{e:3.18a} \quad &&q^{1-\beta}h'_a(q)\nonumber\\
&&\qquad=
\frac{q^{1-\beta}}{(q-1)^{1-\beta}} \int_0^1 \int
_0^{a(q-1)r} \frac{ (-\ln(1-{y}/{(q-1)})
)}{{y}/{(q-1)}} \biggl(1-
\frac{y}{q-1} \biggr)^{q-1}
\\
&&\hspace*{108pt}\qquad\quad{} \times r^\beta \frac{g({y}/{(r(q-1))}) ({y}/{(r(q-1))})^\beta}{y^\beta}
\,dy \,dr.\nonumber
\end{eqnarray}
Hence, again \eqref{e:A}, \eqref{e:A2} and the dominated convergence
theorem yield
%
\begin{equation}
\lim_{q\to\infty} q^{1-\beta}h'_a(q)=
\frac{A\Gamma(1-\beta)}{1+\beta}. \label{e:3.19}
\end{equation}
Here, we use the facts that $(1-\frac{y}{q-1})^{q-1}\le e^{-y}$, $-\ln
(1-z)/z \to1$ as $z\to0$, and also\vspace*{1pt} that $\sup_{z\leq ar<1/2} -\ln
(1-z)/z$ is a finite quantity.
Now
\eqref{e:3.15}, \eqref{e:3.18} and \eqref{e:3.19} jointly imply
\eqref{e:3.11}.

The expression \eqref{e:3.18a} and the bounds just used in deriving
\eqref{e:3.11} also imply that the function $q\mapsto q^{1-\beta
}h'_a(q)$ is bounded on $[2,\infty)$ and, due to the global continuity
of $h_a'$, we conclude that the same function is bounded on $[1,\infty
)$. Together with \eqref{e:3.18} and \eqref{e:3.15}, this proves
\eqref{e:3.12}.
\end{pf}

\section{Integral equations for $N$}
\label{S:Basic}
In this section, we give a representation of the block counting process
$N$ of a given $\Lambda$-coalescent in terms of an integral equation
involving the corresponding Poisson random measure. We also write an
equation for the process $N$ divided by the speed of CDI. Some
preliminary estimates are included at the end.

This construction is our starting point to the proof of the main theorem.
The approach presented here is quite general, and we hope it
to be of independent interest.

In this section and the rest of the paper, we again assume that
$\Lambda
(\{0\})=\Lambda(\{1\})=0$ and that any (and therefore all) of \eqref
{E:cond sch}, \eqref{e:cdi}, \eqref{e:cdi2} hold.

As discussed in Section~\ref{ss:main_tools}, the $\Lambda$-coalescent
can be constructed via a coloring procedure which is based on a Poisson
random measure $\pi$ on $[0,\infty)\times[0,1]$, and an independent
assignment of colors to the blocks. Here, we introduce an enriched
Poisson random measure which contains all the
information on the coloring.
This is a key ingredient in the first important novelty of our
approach---an explicit representation of the martingale which drives the block
counting process $N$.

In order to explain this now,
we will need some additional notation. As usual, let $\N$ denote the
set of natural numbers (without zero).
Let $\mu$ be the law of a sequence of i.i.d. random variables
$X_1,X_2,\ldots$ uniformly distributed on $[0,1]$, that is, $\mu$ is a
probability measure on $[0,1]^\N$, equipped with the product $\sigma
$-algebra generated by the cylinder sets of the form $B_1\times
B_2\times\cdots\times B_n\times[0,1] \times[0,1]\times\cdots,$ $n\in\N$,
$B_i\in\B([0,1])$, $i\in\N$.
The vectors in $[0,1]^\N$ will be denoted in boldface $\x
=(x_1,x_2,\ldots)\in[0,1]^{\N}$. We will usually write $d\x$ instead
of $\mu(d\x)$.

Let $\pi^E$ be a Poisson random measure on $[0,\infty)\times
[0,1]\times
[0,1]^\N$ with intensity measure $ds\frac{\Lambda(dy)}{y^2}\,d\x$.
Observe that such a random measure can be constructed using a Poisson
random measure $\pi$
from (\ref{DPPPpi}) and an independent array of i.i.d. random variables
$(X^i_j)_{i,j\in\N}$, where $X^i_j$ have uniform distribution on
$[0,1]$. Then $\pi^E= \sum_{i\in\N}\delta_{(T_i,Y_i, \bf X^i)}$ is
a Poisson random measure with intensity $ds\frac{\Lambda(dy)}{y^2}\,d\x$.

Moreover, $\pi$ and $\pi^E$ are coupled by the relation
%
\begin{equation}
\pi(\cdot)=\pi^E \bigl(\cdot\times[0,1]^\N \bigr).
\label{e:pi}
\end{equation}
We will henceforth assume that (\ref{e:pi}) holds.
Then we can construct the $\Lambda$ coalescent by the following
procedure: upon arrival of an atom $(t,y,\x)$ of $\pi^E$, the $j$th
block present in the configuration at time $t-$ is colored if and only
if $x_j\leq y$.
Once the colors are assigned, in order to form the configuration at
time $t$, merge all the colored blocks into a single block, and leave
the other (uncolored) blocks intact.

Recall that we assume that the coalescent comes down from infinity, so
$N_r<\infty$ a.s. for any $r>0$. The procedure described above implies that
%
\begin{eqnarray}
\label{e:eqN1} N_t=N_r-\int_{(r,t]\times[0,1]\times[0,1]^\N}
f(N_{s-}, y,\x)\pi^E (ds\,dy\,d\x )
\nonumber
\\[-8pt]
\\[-8pt]
\eqntext{\mbox{for all } 0<r<t,}
\end{eqnarray}
where $f$ is a function which quantifies the decrease in the number of
blocks during a coalescing event:
%
\begin{equation}
\label{e:f} f(k, y, \x)= \Biggl(\sum_{j=1}^k
\mathbf{1}_{\{x_i\leq y\}}-1 \Biggr)\vee0 =\sum_{j=1}^k
\mathbf{1}_{\{x_i\leq y\}}-1 +\mathbf{1}_{\bigcap
_{j=1}^k\{x_j> y\}}.
\end{equation}
Integration with respect to Poisson random measures is well understood;
the reader is referred, for example, to \cite{PZ}.

Recall (\ref{e:Psi}).
One can easily see that
%
\begin{equation}
\label{e:Psi2} \Psi(k)=\int_{[0,1]\times[0,1]^\N} f(k,y,\x)
\frac{\Lambda
(dy)}{y^2}\,d\x.
\end{equation}
Since $\Psi$ is an increasing function and $N$ a decreasing process,
we have
\[
\int_{(r,t]}\Psi(N_{s-})\,ds\le
\Psi(N_{r}) (t-r)\le N_{r}^2 (t-r),
\]
where the last inequality is due to \eqref{e:3.2}.
We know that $EN_{r}^2 <\infty$ [see, e.g., \eqref{e:2.8}] hence,
\[
E \int_{(r,t]\times[0,1]\times[0,1]^\N} f(N_{s-}, y,\x)\,ds
\frac
{\Lambda
(dy)}{y^2}\,d \x<\infty.
\]
This implies that the integral in \eqref{e:eqN1}
belongs to $L^1$ (see, e.g., Theorem~8.23 in~\cite{PZ}).

As the first step toward the proof of Theorem~\ref{thmm:main}, we have
just shown [see \eqref{e:eqN1} and \eqref{e:Psi2}] the following.

%
\begin{lem}
\label{lem:eqN0}
For any $0<r<t$,
%
\begin{equation}\qquad
\label{e:eqN0} N_t=N_r-\int_r^t
\Psi(N_s)\,ds-\int_{(r,t]\times[0,1]\times[0,1]^\N} f(N_{s-},
y, \x){\hat\pi}^E(ds\,dy\,d\x),
\end{equation}
where ${\hat\pi}^E$ denotes the compensated Poisson random measure
%
\begin{equation}
\label{e:hatpiE} {\hat\pi}^E(ds\,dy\,d\x)=\pi^E(ds\,dy\,d\x)-ds
\frac{\Lambda(dy)}{y^2}\,d\x.
\end{equation}
\end{lem}

%
\begin{rem}
The above representation can be done for $N^{(n)}$, the counting
process of the number of blocks of a $\Lambda$-coalescent starting from
$n$ blocks,
even if the $\La$-coalescent does not come down from infinity.
Moreover, a similar representation exists for $\Xi$-coalescents, and
might be useful in similar type of analysis as done here. For
background on this general class of exchangeable coalescents, we refer
the reader to \cite{schweinsberg_xi,ensaios,bertoin}.
\end{rem}

More importantly, we can write a stochastic integral equation for
$\frac
{N_t}{v_t}$.
Indeed, due to \eqref{e:vt} we have
\[
v_t=v_r-\int_r^t
\Psi(v_s)\,ds,\qquad 0<r<t,
\]
thus,
\[
\frac{1}{v_t}=\frac{1}{v_r}+\int_r^t
\frac{\Psi(v_s)}{v_s^2}\,ds
\]
and, therefore, \eqref{e:eqN1} and a simple application of integration
by parts yield

%
\begin{lem}
\label{lem:eqNv0}
For any $0<r<t$,
%
\begin{eqnarray}
\label{e:eqNv0} \frac{N_t}{v_t}&=&\frac{N_r}{v_r}-\int_r^t
\frac{N_s}{v_s} \biggl(\frac
{\Psi(N_s)}{N_s}-\frac{\Psi(v_s)}{v_s} \biggr)\,ds
\nonumber
\\[-8pt]
\\[-8pt]
\nonumber
&&{}-\int_{(r,t]\times[0,1]\times[0,1]^\N} \frac{f(N_{s-}, y,\x
)}{v_s}{\hat\pi}^E (ds\,dy\,d
\x),
\end{eqnarray}
where ${\hat\pi}^E$ is as in \eqref{e:hatpiE}.
\end{lem}

%
\begin{rem}
A predecessor of this result existed in \cite{bbl1,me_xi}, where the
process of main interest was
$\log{N/v^*}$ instead of $N/v$.
The martingale part was not written down explicitly and, therefore,
could not be used in the precise way that it will be used here.
Note that due to \eqref{e:5.4a}, these previous analyses of $\log
{N/v^*}$ as $t\to0$ apply equivalently to $\log{N/v}$.
\end{rem}

It is natural to continue by investigating the integral with respect to
${\hat\pi}^E$.

%
\begin{lem}
\label{lem:Mtilde}
The process $\tilde M=(\tilde M(t))_{t\ge0}$, where
%
\begin{equation}
\label{e:Mtilde} \tilde M(t)=\int_{[0,t]\times[0,1]\times[0,1]^\N} \frac{f(N_{s-},
y,\x
)}{v_s} {
\hat\pi}^E(ds\,dy\,d\x)
\end{equation}
is a well defined, square integrable martingale with quadratic variation
%
\begin{equation}
\label{e:qvariation2} [\tilde M](t)=\int_{[0,t]\times[0,1]\times[0,1]^\N} \biggl(
\frac
{f(N_{s-}, y,\x)}{v_s} \biggr)^2 \pi^E(ds\,dy\,d\x).
\end{equation}
Moreover, for any $p\in(0,2]$, there exists $C(p)>0$, such that for
all $t>0$
%
\begin{equation}
\label{e:4.11} E\sup_{0\le s\le t}\bigl\vert\tilde M(s)
\bigr\vert^p\le C(p) t^{{p}/2}.
\end{equation}
\end{lem}
%

\begin{pf}
Let us first notice that $f(1,\cdot,\cdot)\equiv0$.
Fix $k\in\N$, $k>0$ and $y\in(0,1)$ and let $\xi_{k,y}$ be distributed
as a binomial random variable ${\rm Bin}(k,y)$. Then it is easy to
derive [see also \cite{bbl1}, Lemma~17(iii) and \eqref{e:3.4}--\eqref{e:3.6}]
%
\begin{eqnarray}
\label{e:4.26}
\nonumber
\bck\bck\int_{[0,1]^\N}
f^2(k,y,\x)\,d\x \bck&=&\bck E[ \xi_{k,y} -
\mathbf{1}_{\{{\xi_{k,y}>0}\}}]^2
\\
\bck&=&\bck E( \xi_{k,y})^2 - 2E\xi_{k,y}+ P({
\xi _{k,y}>0})
\\
\bck&=&\bck k(k-1)y^2 - k(k-1)\int_0^y
\int_0^r (1-u)^{k-2}\,du\,dr.\nonumber
\end{eqnarray}
Hence,
%
\begin{eqnarray}
\label{e:qvariation}&& E\int_0^t\int
_{[0,1]\times[0,1]^\N} \biggl(\frac{f(N_{s-}, y,\x)}{v_s} \biggr)^2
\frac{\Lambda(dy)}{y^2} \,ds \,d\x
\nonumber
\\[-8pt]
\\[-8pt]
\nonumber
&&\qquad\le E\int_0^t\int_0^1
\frac{N_{s-}(N_{s-}-1)}{v_s^2}\Lambda (dy)\,ds\le Ct,
\end{eqnarray}
where the last inequality follows from the second moment estimates in
Lem\-ma~\ref{lem:2.2}(v), and the continuity of $v$.

Due to the standard properties of integrals with respect to the
compensated Poisson random measure
(see, e.g., Theorem~8.23 in \cite{PZ}), \eqref{e:qvariation} now
implies that
$\tilde M$ given by \eqref{e:Mtilde} is a well-defined
square integrable martingale with quadratic variation \eqref
{e:qvariation2}. Moreover,
\[
E[\tilde M](t)=\int_{[0,t]\times[0,1]\times[0,1]^\N} E \biggl(\frac
{f(N_{s-}, y,\x)}{v_s}
\biggr)^2\,ds \frac{\Lambda(dy)}{y^2}\,d \x.
\]
Hence, \eqref{e:4.11} for $p=2$ is a consequence of \eqref
{e:qvariation} and the Doob inequality.
The assertion for $0<p<2$ then follows due to Jensen's inequality.
\end{pf}

The bound \eqref{e:4.11} was already implicit in \cite{bbl1}, at least
for $p=2$, where the infinitesimal variance of an analogous martingale
(the one driving the equation for $\log\frac{N_\cdot}{v_\cdot^*}$) was
carefully estimated, even though that martingale was not as explicitly
expressed there as $\tilde M$ is expressed here.

\begin{rem}
\label{R:non Gaussian}
In view of \eqref{e:4.11} for $p=2$
(which becomes an equality asymptotically as $t\to0$),
the fact that both the rate of convergence in Theorem~\ref{thmm:main}
and the law of the limit process depend on rather fine properties of
the driving measure $\Lambda$ may seem surprising.
Without paying consideration to the size of jumps of $N$ at small
times, these inequalities (asymptotic equalities) may suggests Gaussian
type limits for appropriately rescaled $\tilde M$ (and, therefore, for
$N/v-1$). This indeed turns out to be the case in the setting of the
Kingman coalescent (not treated here, check \cite{aldous_survey} for
the nonfunctional CLT in this setting).
However, one quickly realizes that under assumption \eqref{e:A} the
largest jumps of $\tilde M$ (or better, those of $M$) in $[0,\eps t]$
are of order $\eps^{1/(1+\beta)}$. Moreover, if one assumes that
$\La(dy)= \frac{A}{y^\beta} \,dy$ on $[0,1]$ and
denotes by $\Delta_{\eps t}$ the absolute value of the largest jump of
$M$ in $[0,\eps t]$, then it can be easily verified that
$E (\Delta_{\eps t})^2\geq\eps C(\beta,A,t)$,
so the typical bounds on the maximal jump size, sufficient for the
martingale invariance principle to hold [see, e.g., \cite{ethierkurtz}
Chapter~7, Theorem~1.4(b)], are not satisfied here.
Indeed, the Gaussian scaling is not appropriate and, moreover, the
limiting process will have jumps. The paragraph following Remark~\ref
{rem:2.5}(c) gave further intuition regarding the form of the limit.
\end{rem}

Using \eqref{e:eqNv0} and Lemma~\ref{lem:Mtilde}, one can improve on
\eqref{e:2.6} as follows.

%
\begin{lem}
\label{lem:p_moment}
If the $\Lambda$-coalescent comes down from infinity
then for any $p\in(0,2]$ there exists $0<C(p)<\infty$ such that
%
\begin{equation}
\label{e:p_moment} E \sup_{s\le t} \biggl\vert\frac{N_s}{v_s}-1
\biggr\vert^p\le C(p) t^{p/2}.
\end{equation}
\end{lem}

\begin{pf}
Due to Lemma~\ref{lem:3.1}, we know that for any $s>0$, $\frac
{N_{s}}{v_s} (\frac{\Psi(N_s)}{N_s}-\frac{\Psi(v_s)}{v_s}
)$
has the same sign as $\frac{N_s}{v_s}-1$, hence by Lemmas \ref
{lem:eqN0}, \ref{lem:eqNv0}, \ref{lem:Mtilde}
[after subtracting 1 on both sides of (\ref{e:eqNv0})]
and Lemma~\ref{L:bbl1}
we obtain
%
\begin{equation}
\label{e:4.28} \sup_{r\le s\le t}\biggl\vert\frac{N_s}{v_s}-1\biggr\vert\le2
\biggl(\biggl\vert \frac {N_r}{v_r}-1\biggr\vert + \vert\tilde M_r\vert+ \sup
_{r\le s\le t}\vert\tilde M_s\vert \biggr).
\end{equation}
Now \eqref{e:4.11} implies
\[
E \sup_{r\le s\le t}\biggl\vert\frac{N_s}{v_s}-1\biggr\vert^p\le2
\cdot 3^p \biggl(E\biggl\vert\frac{N_r}{v_r}-1\biggr\vert^p + E
\vert\tilde M_r\vert^p + C(p)t^{{p}/2} \biggr).
\]
Letting $r\to0$, and using \eqref{e:2.6} and once again \eqref
{e:4.11}, we obtain \eqref{e:p_moment}.
\end{pf}

\section{Proof of Theorem \texorpdfstring{\protect\ref{thmm:main}}{1.2}}
\label{S:pf_main}
We start this section by giving the scheme of the proof, including an
informal discussion on why Theorem~\ref{thmm:main} should hold.
Our argument is divided into several lemmas, which are proved
separately in the forthcoming subsections.

The first few steps were carried out in Sections \ref{S:preliminary}
and \ref{S:Basic}, while assuming only that the coalescent comes down
from infinity.
Here, as was already done in the final part of Section~\ref{S:Basic},
we specialize further to the case when $\Lambda$ satisfies
assumption~\eqref{e:A}. Recall that \eqref{e:A} implies CDI. Throughout this
section, we assume \eqref{e:A} without much further mention.

The following result is a consequence of Lemmas \ref{lem:eqNv0}, \ref
{lem:Mtilde} and \ref{lem:p_moment}, where assumption \eqref{e:A} makes
passing to the limit $r\searrow0$ possible in the identity \eqref{e:eqNv0}.

\begin{prop}
\label{prop:eqNv}
We have
%
\begin{equation}
\label{e:Nv} \frac{N_t}{v_t} -1=-\int_0^t
\frac{N_s}{v_s} \biggl(\frac{\Psi
(N_s)}{N_s}-\frac{\Psi(v_s)}{v_s} \biggr)\,ds -
\tilde M_t, \qquad t\ge0,
\end{equation}
almost surely, where $\tilde M$ is defined by \eqref{e:Mtilde}.
\end{prop}

\begin{rem}
In the general case [without assuming \eqref{e:A}], one can similarly
obtain a weaker identity,
where the $L^2$ limit
\[
\lim_{r\to0} \int_r^t
\frac{N_s}{v_s} \biggl(\frac{\Psi
(N_s)}{N_s}-\frac
{\Psi(v_s)}{v_s} \biggr)\,ds
\]
exists and replaces the integral from $0$ to $t$ in \eqref{e:Nv}.
At the moment, we do not know whether $s\mapsto\frac{N_s}{v_s}
(\frac{\Psi(N_s)}{N_s}-\frac{\Psi(v_s)}{v_s} )$ is almost surely
Lebesgue integrable on $[0,t]$ in general.
\end{rem}

If $\mathbf{X}=(X_1, X_2, \ldots)$, where $X_i$, $i=1,2,\ldots$ are
i.i.d. random variables uniformly distributed on $[0,1]$, then due to
the form of $f$ [see \eqref{e:f}] and the law of large numbers
it is clear that, for each fixed $y$,
\[
\lim_{k\to\infty}\frac{f(k,y, \mathbf{X})}{k}=y\qquad \mbox{a.s.}
\]
Accounting for \eqref{e:2.5} and $\lim_{t\to0}v_t=\infty$, one would
expect that for small $t$ $\tilde M$ should be close to a martingale
$M=(M(t))_{t\ge0}$ defined by
%
\begin{equation}
M(t)=\int_{[0,t]\times[0,1]}y \hat\pi(ds \,dy), \label{e:M}
\end{equation}
where $\hat\pi$ is the compensated Poisson random measure $\pi$ [see
\eqref{e:pi}], for example,
%
\begin{equation}
\label{e:hatpi} \hat\pi(ds\,dy)=\pi(ds\,dy)-ds\frac{\Lambda(dy)}{y^2}.
\end{equation}
Note that $M$ is a L\'evy process with the L\'evy measure $\frac
{\Lambda
(dy)}{y^2}$.

The above heuristic indeed turns out to be true. More precisely, we
have the following estimate of the difference between $\tilde M$ and $M$:

%
\begin{lem}
\label{lem:diffM}
There exist $t_0>0$ and $0<C<\infty$ such that for all $0<t\le t_0$
%
\begin{equation}
E\sup_{s\le t} (\tilde M_s-M_s
)^2\le C \bigl(t^2\vee t^{{1}/\beta} \bigr).
\label{e:diffM}
\end{equation}
\end{lem}

Concerning the integral on the right-hand side of \eqref{e:Nv}, we have

%
\begin{lem}
\label{lem:drift}
There exist $t_0>0$ and $0<C<\infty$ such that for all $0<t\le t_0$
%
\begin{equation}
\label{e:drift}\qquad E\sup_{u\le t}\biggl \vert\int_0^u
\frac{N_s}{v_s} \biggl(\frac{\Psi
(N_s)}{N_s}-\frac{\Psi(v_s)}{v_s} \biggr)\,ds -\int
_0^u \biggl(\frac
{N_s}{v_s}-1 \biggr)
v_s h'(v_s)\,ds \biggr\vert\le Ct,
\end{equation}
where $h$ is defined by \eqref{e:h}.
\end{lem}

Let us denote by $X$ the process
%
\begin{equation}
X(t)=\frac{N_t}{v_t}-1,\qquad t>0, X(0)=0. \label{e:X}
\end{equation}
Then
\[
X_\varepsilon= \bigl(\varepsilon^{-{1}/{(1+\beta)}}X(\varepsilon t), t\ge0
\bigr)
\]
is the same as the process $X_\varepsilon$ defined in \eqref{e:Xe}.

\textit{Digression-heuristics}. At this point, it is possible to
explain why the limit process of Theorem~\ref{thmm:main} is of the form
as in \eqref{e:Z} (the longer rigorous argument is given below).
From \eqref{e:3.10} and \eqref{e:3.11}, it is not difficult to see that
for $s$ close to zero we have $v_s h'(v_s)\sim\frac{1}s$.
Proposition~\ref{prop:eqNv} and Lemmas \ref{lem:diffM}--\ref{lem:drift}
then jointly give
\[
X(t) \approx-\int_0^t X(s)\frac{1}s
\,ds - M_t.
\]
Making a change of variables in the drift part, we would then have
\[
X_\varepsilon(t)\approx- \int_0^t
X_\varepsilon(s)s^{-1}\,ds - M_{\varepsilon}(t),
\]
where
%
\begin{equation}
\label{e:Me} M_\varepsilon(t)=\varepsilon^{-{1}/{(1+\beta)}}M(\varepsilon t).
\end{equation}

By investigating the Laplace transform of $M_{\varepsilon}$, it is not
difficult to see that it converges in the sense of finite dimensional
distributions to $KL$,
where $L$ is the L\'evy process described in Remark~\ref{rem:2.5}(b)
(this can be verified similarly to Lemma~\ref{lem:fdd} below).
Then it is natural to suspect that, if the limit $Z$ of $X_\varepsilon$
exists, it
should satisfy the equation given in
\eqref{e:eqZ}.
This is indeed the case for the process $Z$ of Theorem~\ref{thmm:main}.

There are a few delicate points in the above reasoning.
We were unable to replace $v_s h'(v_s)$ directly by $\frac{1}s$ and
still get a sufficiently good estimate
(analogous to that of Lemma~\ref{lem:drift})
on the difference between the corresponding integrals.
Furthermore, the convergence of $X_\varepsilon$ has to be proved, and
the passage to the limit under the integral justified.

Our rigorous argument is continued in the following way. Define
%
\begin{equation}
\label{e:Y} Y(t)=\int_{[0,t]}\frac{h(v_t)}{h(v_s)} \,dM(s),\qquad t
\geq0,
\end{equation}
where as usual $h$ is given by \eqref{e:h}, and $M$ by \eqref{e:M}.
We will need the following lemma.

%
\begin{lem}
\label{lem:eqY}
The process $Y$ is the unique solution of the equation
%
\begin{equation}
\label{e:eqY} \,dY(t)=-Y(t)v_th'(v_t)
\,dt+dM(t),\qquad Y(0)=0.
\end{equation}
\end{lem}

Next, we prove that the process $-Y$ is close to $X$.

%
\begin{lem}
\label{lem:approx}
There exist $t_0>0$ and $C>0$ such that
%
\begin{equation}
\label{e:approx} E\sup_{u\le t}\bigl\vert X(u)+Y(u)\bigr\vert\le C \bigl(t
\vee t^{{1}/{(2\beta)
}} \bigr) \qquad\forall t\le t_0.
\end{equation}
\end{lem}

Let $Y_\varepsilon$ denote the following scaled process:
%
\begin{equation}
\label{e:Ye} Y_\varepsilon(t)=\varepsilon^{-{1}/
{(1+\beta)}}Y(\varepsilon t),\qquad
t \geq0.
\end{equation}
Since $1>\frac{1}{1+\beta}$ and $\frac{1}{2\beta}>\frac{1}{1+\beta}$ for
$0<\beta<1$,
Lemma~\ref{lem:approx} implies that\break
$E\sup_{t\le T}\vert X_\varepsilon(t)+Y_\varepsilon(t)\vert\to0$,
for each fixed $T>0$.
In order to prove Theorem~\ref{thmm:main}, it therefore
suffices to show that, as $\varepsilon\to0$, $Y_\varepsilon$ converges
in law to $-Z$ [$Z$ is as defined in \eqref{e:Z}]
with respect to the Skorokhod topology on $D([0,\infty))$, as
$\varepsilon\to0$.

Here we proceed in the standard way: we first derive the convergence of
finite dimensional distributions via the Laplace transform, and then
prove tightness by means of Aldous' tightness criterion.
Let $Z$ be given in \eqref{e:Z}.

%
\begin{lem}
\label{lem:fdd}
As $\varepsilon\to0$, $Y_\varepsilon$ converges to $-Z$
in the sense of finite dimensional distributions.
\end{lem}

%
\begin{lem}
\label{lem:convD}
We have that $Y_\varepsilon\Longrightarrow-Z$ as $\varepsilon\to0$.
\end{lem}

This final lemma, joint with the discussion following the statement of
Lem\-ma~\ref{lem:approx}, completes the proof of Theorem~\ref{thmm:main}.

\subsection{Proof of Proposition \texorpdfstring{\protect\ref{prop:eqNv}}{4.1}}
Let us subtract $1$ on both sides of \eqref{e:eqNv0} and send $r\to0$.
We will show that the integral on the right-hand side of \eqref{e:Nv}
is well defined, and that for any $t>0$ both the left-hand side and the
right-hand side of \eqref{e:eqNv0} with~$1$ subtracted converge in
$L^2$ to the corresponding random variables in \eqref{e:Nv}. This will
imply that for any fixed $t>0$, equation \eqref{e:Nv} is satisfied a.s.
The processes on both sides of \eqref{e:Nv} are right continuous, hence
they are indistinguishable.

Lemma~\ref{lem:Mtilde} [more precisely, \eqref{e:Mtilde} and \eqref
{e:4.11}] implies that the integral with respect to ${\hat\pi}^E$ converges
in $L^2$ to $\tilde M_t$, while Lemma~\ref{lem:2.2} part (v) implies
that $\frac{N_r}{v_r}-1$ converges to $0$ in $L^2$.
Therefore, the remaining term on the right-hand side of \eqref{e:eqNv0}
must also converge in $L^2$.
Moreover, it is not hard to see that the integral
\[
\int_0^t \frac{N_s}{v_s} \biggl(
\frac{\Psi(N_s)}{N_s}-\frac{\Psi
(v_s)}{v_s} \biggr)\,ds = \int_0^t
\frac{N_s}{v_s} \bigl(h(N_s)-h(v_s) \bigr) \,ds
\]
is well defined a.s. as a Lebesgue integral.
Indeed, the derivative of $h$ is nonnegative due to Lemma~\ref
{lem:3.1} part (iv).
We will repeatedly use assumption \eqref{e:A} in the rest of the argument.
Observe that
\eqref{e:3.15}--\eqref{e:3.17} imply that $h'$ is decreasing.
Hence, if $N_s\le v_s$, then
\begin{eqnarray*}
\frac{N_s}{v_s} \bigl\vert h(N_s)-h(v_s)\bigr\vert &\le&
{N_s}h'(N_s)\biggl \vert\frac{N_s}{v_s}-1
\biggr\vert
\\
&\le& C \biggl(\frac{1}{s}\vee1 \biggr) \biggl\vert\frac{N_s}{v_s}-1\biggr\vert,
\end{eqnarray*}
where the last inequality follows from \eqref{e:3.12}, the fact that
$N_s^\beta\le v_s^\beta$ and \eqref{e:3.10a}.

If $N_s>v_s$, then again by \eqref{e:3.10a} and \eqref{e:3.12}
\begin{eqnarray*}
\frac{N_s}{v_s} \bigl\vert h(N_s)-h(v_s)\bigr\vert &\le&
N_s h'(v_s)\biggl \vert\frac{N_s}{v_s}-1
\biggr\vert
\\
&\le& C \biggl(\frac{1}{s}\vee1 \biggr) \frac{N_s}{v_s} \biggl\vert
\frac
{N_s}{v_s}-1\biggr\vert.
\end{eqnarray*}
The Cauchy--Schwarz inequality, Lemma~\ref{lem:p_moment} and \eqref
{e:2.8} now imply that
\begin{eqnarray*}
E \biggl(\int_0^t \frac{N_s}{v_s}\biggl\vert
\frac{\Psi(N_s)}{N_s}-\frac
{\Psi (v_s)}{v_s}\biggr\vert \,ds \biggr)&\leq& C E \int
_0^t \biggl(\frac{1}{s} \vee1 \biggr)
\biggl(1+\frac
{N_s}{v_s} \biggr) \biggl\llvert \frac{N_s}{v_s}-1 \biggr
\rrvert \,ds
\\
&\leq& C_1 \int_0^t \biggl(
\frac{1}{s} \vee1 \biggr) \sqrt{s} \,ds<\infty.
\end{eqnarray*}
Letting $r\to0$ in \eqref{e:eqNv0}, we obtain \eqref{e:Nv}.

\subsection{Proof of Lemma \texorpdfstring{\protect\ref{lem:diffM}}{4.3}}
Recalling\vspace*{1pt} the forms of $M$ and $\tilde M$ [see \eqref{e:M} and~\eqref
{e:Mtilde}] as well as \eqref{e:pi},
observe that $\tilde M-M$ is a square integrable martingale with
quadratic variation process
\[
[\tilde M-M](t)=\int_{[0,t]\times[0,1]\times[0,1]^\N} \biggl(\frac
{f(N_{s-},y,\x)}{v_s}-y
\biggr)^2 \pi^E(ds\,dy\,d\x).
\]
Thus, we have
\[
E[\tilde M - M](t)\le2EI_1(t)+2EI_2(t),
\]
where
%
\begin{equation}
\label{e:4.30} I_1(t)=\int_{[0,t]\times[0,1]\times[0,1]^\N} \biggl(
\frac
{f(N_{s-},y,\x
)-N_{s-}y}{v_s} \biggr)^2 \,ds\frac{\Lambda(dy)}{y^2}\,d\x
\end{equation}
and
%
\begin{equation}
\label{e:4.31} I_2(t)=\int_{[0,t]\times[0,1]} \biggl(
\frac{N_s}{v_s}-1 \biggr)^2 \,ds\Lambda(dy)= \int
_0^t \biggl(\frac{N_s}{v_s}-1
\biggr)^2 \,ds.
\end{equation}
By Doob's inequality, it therefore suffices to show
%
\begin{equation}
\label{e:4.32} EI_i(t)\le C \bigl(t^2\vee
t^{{1}/\beta} \bigr),\qquad i=1,2.
\end{equation}
Estimate \eqref{e:4.32} for $I_2$ is immediate by Lemma~\ref{lem:p_moment}.
Arguing \eqref{e:4.32} for $I_1$ is a bit more involved. Let us denote
%
\begin{equation}
\label{e:4.33} J(k)=\int_0^1 \int
_{[0,1]^\N} \bigl(f(k,y,\x)-ky \bigr)^2 \,d\x
\frac
{\Lambda(dy)}{y^2},\qquad k\in\N,
\end{equation}
so that
\[
I_1(t) = \int_{[0,t]} \frac{J(N_{s-})}{v_s^2}\,ds.
\]
By \eqref{e:f}, \eqref{e:4.26} and the following, easy to check identity
\[
\int_{[0,1]^\N}f(k,y,\x)\,d\x=ky-k\int_0^y
(1-r)^{k-1}\,dr,
\]
we have
\[
J(k)\le2k^2 \int_0^1 \int
_0^y (1-r)^{k-1} \cdot y \,dr
\frac
{\Lambda
(dy)}{y^2}.
\]
Taking $a$ which satisfies \eqref{e:A2}, and applying $1-r\le e^{-r}$
we write
%
\begin{equation}
\label{e:4.34} J(k)\le2e \bigl(J_a(k)+\tilde J_a(k)
\bigr),
\end{equation}
where
\[
J_a(k)=k^2\int_0^a
\int_0^y e^{-kr}\,dr
\frac{\Lambda(dy)}{y},\qquad \tilde J_a(k)=k^2\int
_a^1 \int_0^y
e^{-kr}\,dr \frac{\Lambda(dy)}{y}.
\]
By \eqref{e:A2} and the natural substitutions ($r'=r/y$, followed by
$y'=kr'y$, and afterward $r', y'$ renamed to $r, y$, resp.) we have
\[
J_a(k)\le C k^{1+\beta} \int_0^1
\int_0^{akr} e^{-y}y^{-\beta}
r^{\beta
-1}\,dy \,dr \le C_1 k^{1+\beta}.
\]
The term $\tilde J_a$ can be easily bounded as follows:
\[
\tilde J_a(k) \le\frac{k}{a}.
\]
Recalling \eqref{e:4.34}, we therefore have
$J(k)\le C k^{1+\beta}$ for some $C<\iy$.
Together with \eqref{e:4.33}, \eqref{e:4.30}, \eqref{e:2.8} and
\eqref
{e:3.10a}, this now implies that
(for $t_0<1/2$ we use $1\vee1/s= 1/s$, $\forall s<t_0$)
\begin{eqnarray*}
EI_1(t)&=& E\int_0^t
J(N_{s-}) \frac{1}{v_s^2} \,ds \le C E\int_0^t
\biggl(\frac{N_s}{v_s} \biggr)^{1+\beta}v_s^{\beta-1}
\,ds
\\
&\le& C_1\int_0^t \biggl(
\frac{1}{s^{1/\be}} \biggr)^{\be-1} \,ds = C_2t^{1/\be},
\end{eqnarray*}
which proves \eqref{e:4.32} for $i=1$, and completes the argument.

\subsection{Proof of Lemma \texorpdfstring{\protect\ref{lem:drift}}{4.4}}
Let $h$ be defined by \eqref{e:h} and let $h_a$ and $\tilde h_a$ be as
in \eqref{e:3.15a}--\eqref{e:3.15b}, with $0<a<\frac{1}2$ satisfying
\eqref{e:A2}. Using the
easy estimate $\tilde h_a(q)\le a^{-2}$ together with \eqref{e:2.8},
we have
\[
E\int_0^t \frac{N_s}{v_s}\bigl\vert\tilde
h_a(N_s)- \tilde h_a(v_s)
\bigr\vert \,ds\le Ct.
\]
Moreover, by \eqref{e:3.18}, Lemma~\ref{lem:p_moment} and \eqref
{e:3.10a} we obtain
\[
E\int_0^t \biggl\vert\frac{N_s}{v_s}-1\biggr\vert
v_s\tilde h_a'(v_s)\,ds\le
Ct^{{1}/\beta+ {3}/2}.
\]
Hence, to prove the lemma, it suffices to show \eqref{e:drift} with $h$
replaced by $h_a$.
Using the Taylor expansion formula, we write
%
\begin{equation}
\label{e:4.37} \frac{N_s}{v_s} \bigl(h_a(N_s)-h_a(v_s)
\bigr)=I_1(s)+I_2(s),
\end{equation}
where
\[
I_1(s)=\frac{N_s}{v_s}\frac{N_s-v_s}{v_s} v_sh_a'(v_s),\qquad
I_2(s)=\frac
{N_s}{v_s}\int_{v_s}^{N_s}
\int_{v_s}^z h_a''(w)
\,dw\,dz.
\]
We shall prove that $I_1$ is the main term, uniformly close to
$(N_\cdot-v_\cdot)h_a'(v_\cdot)$, and that $I_2$ is a negligible
error term.
First note that
by Lemma~\ref{lem:p_moment}, \eqref{e:3.12} (recall that $h'_a\le h'$)
and \eqref{e:3.10a} one can easily see that
%
\begin{equation}
\label{e:basich} E \biggl\llvert \biggl(\frac{N_s}{v_s}-1 \biggr)
(N_s - v_s) h_a'(v_s)
\biggr\rrvert \leq C E \biggl(\frac{N_s}{v_s}-1 \biggr)^2
v_s^\beta= O(1),
\end{equation}
and, therefore,
%
\begin{equation}
\label{e:4.40} E\int_0^t\bigl\vert
I_1(s)-(N_s-v_s)h'_a(v_s)
\bigr\vert \,ds \le Ct.
\end{equation}

Our approach for $I_2$ is to show a similar bound
%
\begin{equation}
\label{e:I2} \bigl\vert I_2(s)\bigr\vert\le C \biggl(\frac{N_s-v_s}{v_s}
\biggr)^2 v_s^\beta,
\end{equation}
and then again use \eqref{e:basich} to bound $\int_0^t |I_2(s)| \,ds$.
First, note that from differentiating in
\eqref{e:3.16} it follows that $h''_a$ is negative and increasing (its
absolute value is decreasing).
Moreover, since $a<\frac{1}2$, and since $|\log(1-r)|\leq2r$ and
$(1-r)^{q-1} \leq2 e^{-rq}$ for $r\leq1/2$, one can easily
derive from \eqref{e:A2} that
%
\begin{equation}
\bigl\vert h''_a(q)\bigr\vert\le C \int
_0^a\int_0^y
r^2 e^{-rq} y^{-2-\beta
}\,dr \,dy =O \bigl(
q^{\beta-2} \bigr). \label{e:4.41}
\end{equation}
Thus, if $\frac{1}2 v_s\le N_s \le2v_s$, then
\[
\bigl|h_a''(w)\bigr| \le\biggl\vert h_a''
\biggl(\frac{1}2 v_s \biggr)\biggr\vert =O \bigl(
v_s^{\be
-2} \bigr)
\]
and
$\vert I_2(s)\vert= \frac{N_s}{v_s}(N_s-v_s)^2O(v_s^{\beta-2})$.
Since $N_s/v_s \le2$, we conclude that \eqref{e:I2} holds in this case.

If $v_s>2N_s$ then note that
\begin{eqnarray*}
\int_{v_s}^{N_s} \int_{v_s}^{z}
w^{\be-2} \,dw \,dz &=& \int_{N_s}^{v_s}\int
_{z}^{v_s} w^{\be-2} \,dw \,dz
\\
&\leq&\frac{1}{1-\be} \int_{N_s}^{v_s}
z^{\be-1} \,dz
\\
&\leq&\frac{1}{1-\be}(v_s-N_s)N_s^{\beta-1}.
\end{eqnarray*}
Hence, by \eqref{e:4.41} and the definition of $I_2$
\[
\bigl\vert I_2(s)\bigr\vert\le C \biggl(\frac{v_s-N_s}{v_s} \biggr)
N_s^\beta.
\]
We also have $N_s^\beta\le v_s^\beta$ and $1 <2 \frac{v_s-N_s}{v_s}$,
so \eqref{e:I2} follows.

If $2v_s<N_s$, then
\begin{eqnarray*}
\frac{N_s}{v_s}\int_{v_s}^{N_s}\int
_{v_s}^z w^{\beta-2}\,dw \,dz &\le& C
\frac{N_s}{v_s} (N_s-v_s)v_s^{\beta-1}
\\
&\le& C \biggl(\frac{N_s}{v_s}-1 \biggr)^2 v_s^\beta+C
\biggl(\frac
{N_s}{v_s}-1 \biggr) v_s^\beta.
\end{eqnarray*}
Together with \eqref{e:4.41} and the definition of $I_2(s)$ this again implies
\eqref{e:I2}, since for $2v_s<N_s$ we have $1<\frac{N_s}{v_s}-1<
(\frac{N_s}{v_s}-1 )^2$.

This gives \eqref{e:I2}, and due to the final estimate in \eqref
{e:basich} we get\break
$ E\int_0^t \vert I_2(s)\vert \,ds\le Ct$, which combined with \eqref{e:4.40}
yields\eqref{e:drift} for $h_a$.
As already argued, this completes the proof of the lemma.

\subsection{Proof of Lemma \texorpdfstring{\protect\ref{lem:eqY}}{4.5}}
Let us first observe that the function $u\mapsto h(v_u)$
defined in \eqref{e:h} is positive on $(0,\infty)$ and strictly
decreasing, since
$h$ is positive and strictly increasing and $v$ is strictly decreasing
(see Lemmas \ref{lem:3.1} and \ref{lem:2.2}). Moreover, by \eqref
{e:3.9} and \eqref{e:3.10}, we have that
%
\begin{equation}
\label{e:4.44} \lim_{u\to0}{uh(v_u)}=
\frac{1}{\beta},
\end{equation}
so, there exists $t_0$ such that
%
\begin{equation}
\label{e:4.45} \frac{\beta}{2} u\le\frac{1}{h(v_u)}\le2\beta u,\qquad 0<u\le
t_0.
\end{equation}
Hence, the process $Y$ from \eqref{e:Y} is well defined. Moreover,
%
\begin{equation}
\label{e:4.46} E \bigl(Y(t) \bigr)^2= \bigl(h(v_t)
\bigr)^2 \int_0^t\int
_0^1 \biggl(\frac
{y}{h(v_u)}
\biggr)^2 \frac{\Lambda(dy)}{y^2}\le t,
\end{equation}
since $h(v_t)\le h(v_u)$ for $u\le t$.

The function $u\mapsto h(v_u)$ is clearly continuous and of finite
variation on any interval $[r,t]$, $0<r<t$.
We apply integration by parts, which in this case is simply
$fg=\int f \,dg +\int g \,df$ with $f(\cdot)=h(v_\cdot)$ and $g(\cdot
)=\int_0^\cdot\frac{1}{h(v_s)}\,dM_s$
(note that the other terms which normally appear in this formula are
equal to 0, due to just mentioned continuity and finite variation properties).
Using the fact that
$ \frac{v_s'}{h(v_s)}=-v_s$,
[cf. \eqref{e:h} and \eqref{e:2.3a}], we get for $0<r<t$
%
\begin{equation}
Y_t=Y_r-\int_r^t
Y_sv_sh'(v_s)\,ds
+M_t-M_r. \label{e:4.47}
\end{equation}
We now let $r\to0$ and observe that $M_r\to0$ a.s. and in $L^2$,
since $E[M](r)=\int_0^r\int_0^1y^2 \frac{\Lambda(dy)}{y^2}=r$, and
$Y_r\to0$ in $L^2$ by \eqref{e:4.46}. To deal with the remaining term
in \eqref{e:4.47}, we note that by \eqref{e:3.12} and \eqref{e:3.10a}
we have
\[
0\le v_s h'(v_s)\le C
\bigl(s^{-1}\vee1 \bigr).
\]
Hence, by \eqref{e:4.46} and Jensen's inequality
\[
E\int_0^r\bigl\vert Y_s
v_s h'(v_s)\bigr\vert \,ds\le C \int
_0^r \sqrt{s} \biggl( \frac{1}{s} \vee1
\biggr) \,ds \le C \bigl(\sqrt{r}\vee r^{3/2} \bigr),
\]
converges to $0$ as $r\to0$. After sending $r\to0$ in \eqref{e:4.47},
one concludes that $Y$ given by \eqref{e:Y} satisfies equation \eqref{e:eqY}.

Showing uniqueness is easier. Indeed, if $Y_1$ and $Y_2$ are two
solutions of \eqref{e:eqY}, then
\[
Y_1(t)-Y_2(t)=-\int_0^t
\bigl(Y_1(s)-Y_2(s) \bigr)v_sh'(v_s)
\,ds.
\]
Since $v_sh'(v_s)$ is positive [see Lemma~\ref{lem:3.1}(iv)], an
application of Lemma~\ref{L:bbl1} implies $Y_1-Y_2\equiv0$.

\subsection{Proof of Lemma \texorpdfstring{\protect\ref{lem:approx}}{4.6}}

Recall \eqref{e:X}.
Due to Proposition~\ref{prop:eqNv} and Lemmas~\ref{lem:diffM}, \ref
{lem:drift} and \ref{lem:eqY}, we obtain
\[
X(t)+Y(t)=-\int_0^t \bigl(X(s)+Y(s)
\bigr)v_s h'(v_s)\,ds +R(t),
\]
where $R$ is a process such that for $0\le t\le t_0$
\[
E\sup_{s\le t}\bigl\vert R(s)\bigr\vert\le C \bigl(t\vee
t^{{1}/ {(2\beta) }} \bigr).
\]
Since $v_s h'(v_s)$ is positive, another application of Lemma~\ref
{L:bbl1} completes the proof.

\subsection{Proof of Lemma \texorpdfstring{\protect\ref{lem:fdd}}{4.7}}

The argument relies on convergence of the Laplace transform for
positive arguments.
Fix $n\in\N$ and $z_j\ge0$, $t_j>0$, $j=1,2,\ldots, n$ and denote
%
\begin{equation}
\label{e:F} F(u)=\sum_{j=1}^n
z_j \frac{u}{t_j}\ind_{[0,t_j]}(u).
\end{equation}
We will show that
%
\begin{eqnarray}
\label{e:4.48}&& \lim_{\varepsilon\to0} E\exp \Biggl\{- \sum
_{j=1}^n z_jY_\varepsilon
(t_j) \Biggr\}
\nonumber
\\[-8pt]
\\[-8pt]
\nonumber
&&\qquad= \exp \biggl\{ A\int_0^\infty
\bigl(e^{-y}-1+y \bigr)\frac{1}{y^{2+\beta}}\,dy \int_0^\infty
\bigl(F(u) \bigr)^{1+\beta}\,du \biggr\}.
\end{eqnarray}
Due to Propositions 3.4.1 and 1.2.12 and (3.4.4) in \cite{ST}, the
right-hand side is precisely
$E\exp{  \{- \sum_{j=1}^n z_j(-Z(t_j)) \}}$, where $Z$ is
defined in \eqref{e:Z}.
On the other hand, it is well known that since $-Z$ is a $(1+\beta
)$-stable process totally skewed to the right, the convergence of
Laplace transforms for all positive $z_j$ implies the convergence in
law of $(Y_\varepsilon(t_1),\ldots,Y_\varepsilon(t_n))$ to
$(-Z(t_1),\ldots,-Z(t_n))$ (see, e.g., \cite{Iscoe}, proofs of Theorems
5.4 and 5.6).
Thus, the lemma will be proved once we show~\eqref{e:4.48}.

By \eqref{e:Y} and \eqref{e:Ye}, we have
\begin{eqnarray*}
\sum_{j=1}^n z_jY_\varepsilon(t_j)
&=&\varepsilon^{-{1}/{(1+\beta)
}}\int_0^\infty\int
_0^1 \Biggl(\sum_{j=1}^n
z_j \ind_{[0,\varepsilon t_j]}(u) \frac
{h(v_{\varepsilon t_j})}{h(v_u)} \Biggr)y \hat
\pi(du\,dy)
\\
&=&\varepsilon^{-{1}/{(1+\beta)}}\int_0^\infty\int
_0^1 F_\varepsilon \biggl(
\frac{u}\varepsilon \biggr) y \hat\pi(du\,dy),
\end{eqnarray*}
where
%
\begin{equation}
\label{e:Fe} F_\varepsilon(u)=\sum_{j=1}^n
z_j \frac{h(v_{\varepsilon
t_j})}{h(v_{\varepsilon u})}\ind_{[0,t_j]}(u).
\end{equation}
Thus, by the usual properties of a Poisson random measure, we have
%
\begin{equation}
\label{e:4.49} E\exp \Biggl\{ -\sum_{j=1}^n
z_jY_\varepsilon(t_j) \Biggr\} =e^{I(\varepsilon)},
\end{equation}
where
%
\begin{equation}\qquad
\label{e:4.50} I(\varepsilon)=\int_0^\infty\int
_0^1 \biggl(e^{-\varepsilon
^{-{1}/{(1+\beta)}}F_\varepsilon({u}/\varepsilon)y} -1+
\varepsilon^{-{1}/{(1+\beta)}}F_\varepsilon \biggl(\frac
{u}{\varepsilon} \biggr)y
\biggr) \frac{\Lambda(dy)}{y^2} \,du.\hspace*{-6pt}
\end{equation}
As before, let $0<a<\frac{1}2$ be such that \eqref{e:A2} holds and write
%
\begin{equation}
\label{e:4.51} I(\varepsilon)=I_a(\varepsilon) +\tilde
I_a(\varepsilon),
\end{equation}
where
%
\begin{equation}
\label{e:4.52} I_a(\varepsilon)=\int_0^\infty
\int_0^a \cdots\quad \mbox{and}\quad \tilde
I_a(\varepsilon)=\int_0^\infty\int
_a^1 \cdots,
\end{equation}
and the $\cdots$ above denotes the expression under the integral in
\eqref{e:4.50}.
Let us initially consider the term $\tilde I_a$. We have
\begin{eqnarray*}
0&\le&\tilde I_a(\varepsilon)\bck\le\bck\int_0^\infty
\int_a^1 \varepsilon^{-{1}/{(1+\beta)}}
F_\varepsilon \biggl(\frac{u}\varepsilon \biggr)\frac{\Lambda(dy)}{y}
\,du
\\
\bck&\le&\bck\frac{1}a \varepsilon^{1-{1}/{(1+\beta)}}\int
_0^\infty F_\varepsilon(u)\,du.
\end{eqnarray*}
Recall \eqref{e:Fe} and note that
$h(v_{\varepsilon t})\le h(v_{\varepsilon u})$ for $u\le t$, as
explained in the proof of Lemma~\ref{lem:eqY}. Thus, $\sup_{\varepsilon
>0}\int_0^\infty F_\varepsilon(u)\,du<\infty$ and
it follows that
%
\begin{equation}
\label{e:4.54} \lim_{\varepsilon\to0} \tilde I_a(
\varepsilon)=0.
\end{equation}
In the analysis of $I_a(\eps)$, we make a change of variables $y=z\eps
^{1/(1+\be)}$ and $r=\frac{u}{\varepsilon}$
(then rename $z$ to be $y$ and $r$ to be $u$) and use assumption \eqref
{e:A} to get
%
\begin{eqnarray}\label{e:4.55}
&&I_a(\varepsilon)=\int_0^\infty\int
_0^{a\varepsilon^{-
{1}/{(1+\beta)}}} \bigl(e^{-F_\varepsilon(u)y}-1+F_\varepsilon(u)y
\bigr)
\nonumber
\\[-8pt]
\\[-8pt]
\nonumber
&&\hspace*{104pt}{}\times\frac{g(y\varepsilon^{{1}/{(1+\beta)}})  ( y\varepsilon
^{{1}/{(1+\beta)}} )^\beta} {
y^{2+\beta}} \,dy \,du.
\end{eqnarray}
By \eqref{e:4.44},
we have
\[
\lim_{\varepsilon\to0}\frac{h(v_{\varepsilon t})}{h(v_{\varepsilon
u})}=\frac{u}{t},
\]
so from \eqref{e:Fe} we see that
$F_\varepsilon$ converges pointwise to $F$ defined in \eqref{e:F}.
Moreover, note that
\[
0\le e^{-F_\varepsilon(u)y}-1+F_\varepsilon(u)y\le F^2_\varepsilon
(u)y^2\le \Biggl(\sum_{j=1}^n
z_j\ind_{[0,t_j]}(u) \Biggr)^2 y^2.
\]
Hence, by \eqref{e:4.55}, \eqref{e:A}, \eqref{e:A2} and the dominated
convergence theorem, it follows that
%
\begin{eqnarray}
\label{e:4.56} \lim_{\varepsilon\to0}I_a(\varepsilon) &= &A
\int_0^\infty\int_0^\infty
\bigl(e^{-F(u)y}-1+F(u)y \bigr)\frac
{1}{y^{2+\beta}}\,dy \,du
\nonumber
\\[-8pt]
\\[-8pt]
\nonumber
&= &A\int_0^\infty \bigl(F(u)
\bigr)^{1+\beta}\,du\int_0^\infty
\bigl(e^{-y}-1+y \bigr)\frac{1}{y^{2+\beta}}\,dy,
\end{eqnarray}
where we apply the substitution $z=F(u)y$ and then rename $z$ as $y$.
Now \eqref{e:4.49}--\eqref{e:4.52}, \eqref{e:4.54} and \eqref{e:4.56}
together imply that \eqref{e:4.48} holds and the proof is complete.

\subsection{Proof of Lemma \texorpdfstring{\protect\ref{lem:convD}}{4.8}}
First observe that by \eqref{e:4.44} and \eqref{e:4.45} the function
$f_\varepsilon$ defined by $f_\varepsilon(0)=\frac{1}\beta$ and
$f_\varepsilon(t)=\varepsilon t h(v_{\varepsilon t})$ for $t>0$
is continuous for any $\varepsilon>0$.
Furthermore, as $\varepsilon\to0 $, the family $(f_\varepsilon
)_{\varepsilon>0}$ converge uniformly on bounded intervals to a
constant function $\frac{1}\beta$.
Hence, to prove the lemma, it suffices to show that the family of processes
$(\tilde Y_\varepsilon)_\varepsilon$ defined by
%
\begin{equation}
\tilde Y_\varepsilon(t)=t^{-1} \beta^{-1}
\varepsilon^{-1-{1}/{(1+\beta)}}\int_0^{\varepsilon t}
\frac{1}{h(v_u)}\,dM_u \label{e:tildeYe}
\end{equation}
converges in law in $D([0,\infty))$ to $-Z$, as $\varepsilon\to0$.

We will split the proof into several steps. In the first step, with the
help of
Aldous' tightness criterion, we show that the family of processes
$(t\tilde Y_\varepsilon(t))_{t\ge0}$ converges in law in $D([0,\infty
))$ to $(-tZ(t))_{t\ge0}$.
From this, we need to infer the convergence $\tilde Y_\varepsilon
\Rightarrow-Z$. However, the latter step is not immediate, since the
function $t\mapsto\frac{1}t$ cannot be extended to a continuous
function on $[0,\infty)$. We will overcome this problem by taking
suitable approximations.

\textit{Step} 1.
We prove that the family of processes $(U_{\varepsilon})_{\varepsilon
>0}$ defined by
%
\begin{equation}
\label{e:U} U_{\varepsilon}(t)=t\tilde Y_\varepsilon(t),\qquad t\ge0,
\varepsilon>0,
\end{equation}
converges to $(-tZ(t))_{t\ge0}$ in law in $D([0,\infty))$.
It is clearly enough to show this convergence when restricted to an
arbitrary but fixed sequence $\varepsilon_n\searrow0$.

The convergence of finite dimensional distributions follows from \eqref
{e:4.44} and Lemma~\ref{lem:fdd}.
To prove tightness of the family $(U_{\varepsilon})_{\varepsilon>0}$,
we will apply
the well-known Aldous criterion (see, e.g., \cite{Billingsley} Theorem~16.10).
More precisely, we will prove:

(i) For any $M>0$,
%
\begin{equation}
\lim_{r\to\infty}\limsup_{n\to\infty}\P \Bigl(\sup
_{t\in[0,M]}\bigl\vert U_{\varepsilon_n}(t)\bigr\vert\ge r \Bigr)=0,
\label{e:Aldous1}
\end{equation}

(ii) For any $\rho,\eta, M>0$, there exist $\delta_0, n_0$ such that if
$\delta\le\delta_0$, $n\ge n_0$ and $\tau$ is a stopping time with
respect to the filtration generated by $U_{\varepsilon_n}$, taking
finite number of values, and
such that $\P(\tau\le M)=1$, then
%
\begin{equation}
\P \bigl(\bigl\vert U_{\varepsilon_n}(\tau+\delta)-U_{\varepsilon_n}(\tau )
\bigr\vert\ge
\rho \bigr)\le\eta. \label{e:Aldous2}
\end{equation}

To prove (i) and (ii) above, we will need an estimate on the moments of
increments of $U_{\varepsilon}$.
We write
%
\begin{equation}
\label{e:U12} U_\varepsilon=\frac{1}\beta \bigl(U^{(1)}_{\varepsilon}+U_\varepsilon
^{(2)} \bigr),
\end{equation}
where
\begin{eqnarray*}
U_{\varepsilon}^{(1)}(r)&=&\varepsilon^{-{(2+\beta)}/{(1+\beta)}} \int
_{[0,\varepsilon r]\times[0,\varepsilon^{{1}/{(1+\beta)}}]} \frac{1}{h(v_u)}y\hat\pi(du\,dy),
\\
U_{\varepsilon}^{(2)}(r)&=&\varepsilon^{-{(2+\beta)}/{(1+\beta)}} \int
_{[0,\varepsilon r]\times(\varepsilon^{{1}/{(1+\beta)}},1]} \frac{1}{h(v_u)}y\hat\pi(du\,dy).
\end{eqnarray*}
Note that $U^{(1)}_\varepsilon$ (resp., $U^{(2)}_\varepsilon$) is the process
which captures the ``small'' (resp., ``large'')
jumps of $U_\varepsilon$.

Using standard properties of integrals with respect to a compensated
Poisson random measure (see, e.g., \cite{PZ}, Theorem~8.23), we have
\[
E\bigl\vert U_\varepsilon^{(2)}(t) - U_\varepsilon^{(2)}(s)
\bigr\vert^p\le C\varepsilon ^{-{p(2+\beta)}/{(1+\beta)}}\int_{\varepsilon s}^{\varepsilon
t}
\int_{\varepsilon^{{1}/{(1+\beta)}}}^1 \frac{y^p}{(h(v_u))^p}
\frac
{\Lambda
(dy)}{y^2}\,du.
\]
Let $0<s<t<T$ and $1<p<1+\beta$ and suppose that $\varepsilon\le
{a^{1+\beta}} \wedge\frac{t_0} T$, where $a$ is as in \eqref{e:A2}
and $t_0$ as in \eqref{e:4.45}.
By \eqref{e:A2} and \eqref{e:4.45}, we obtain
%
\begin{eqnarray}\label{e:4.57}
&&E\bigl\vert U_\varepsilon^{(2)}(t) - U_\varepsilon^{(2)}(s)
\bigr\vert ^p
\nonumber
\\
&&\qquad\le C \varepsilon^{-{p(2+\beta)}/{(1+\beta)}}\int_{\varepsilon
s}^{\varepsilon t}u^p
\biggl(\int_{\varepsilon^{{1}/{(1+\beta)}}}^{a} y^{p-2-\beta}\,dy +\int
_{a}^{1} y^{p-2}\Lambda(dy) \biggr)\,du
\\
&&\qquad\le C_1(p) T^p (t-s),\nonumber
\end{eqnarray}
since $\eps^{p+1}\ll\eps^{p+1} \eps^{{(p-1-\be)}/{(1+\be)}}=\eps
^{{p(2+\be)}/{(1+\be)}}$ cancels the power of $\eps$ in front of the
integral, and since
$\int_{a}^{1} y^{p-2}\Lambda(dy)$ is a constant quantity.

Via similar arguments applied to $U^{(1)}$, we get
\[
E\bigl\vert U_\varepsilon^{(1)}(t) - U_\varepsilon^{(1)}(s)
\bigr\vert ^2=\varepsilon ^{-{2(2+\beta)}/{(1+\beta)}}\int_{\varepsilon s}^{\varepsilon
t}
\int_0^{\varepsilon^{{1}/{(1+\beta)}}} \frac{1}{(h(v_u))^2}\Lambda(dy)\,du,
\]
and, since $3+\frac{1-\be}{1+\be}=\frac{2(2+\be)}{1+\be}$, again
\eqref
{e:A2} and \eqref{e:4.45} yield
%
\begin{eqnarray}
\label{e:4.58} E\bigl\vert U_\varepsilon^{(1)}(t) -
U_\varepsilon^{(1)}(s) \bigr\vert^2& \le& C
\varepsilon^{-{2(2+\beta)}/{(1+\beta)}}\int_{\varepsilon
s}^{\varepsilon t} \int
_0^{\varepsilon^{{1}/{(1+\beta)}}} \frac
{u^2}{y^\beta}\,dy \,du
\nonumber
\\[-8pt]
\\[-8pt]
\nonumber
&\le& C_2 T^2 (t-s).
\end{eqnarray}
Now \eqref{e:U12}--\eqref{e:4.58} and Jensen's inequality imply that for
$0<s<t<T$ and $1<p<1+\beta$,
$\varepsilon\le{a^{1+\beta}} \wedge\frac{t_0} T$ we have
%
\begin{equation}
E\bigl\vert U_\varepsilon(t)-U_\varepsilon(s)\bigr\vert^p\le C(p)
T^{p} \bigl( \vert t-s\vert^{{p}/ 2} \vee\vert t-s\vert
\bigr). \label{e:4.59}
\end{equation}

Applying the Doob maximal inequality to the martingale $U_\varepsilon$,
we conclude
\[
\P \Bigl(\sup_{t\in[0,M]}\bigl\vert U_\varepsilon(t)\bigr\vert>r \Bigr)
\le \biggl(\frac{p}{p-1} \biggr)^p \frac{E\vert U_\varepsilon(M)\vert^p}{r^p}.
\]
Hence, \eqref{e:4.59} implies \eqref{e:Aldous1}.

Estimate \eqref{e:4.59} and the
Markov property (since $\tau$ takes only finitely many values, we do
not need the strong Markov property)
of $U_\varepsilon$ imply that
if $\tau$ is a stopping time with respect to the filtration of
$U_\varepsilon$ taking finite number of values and such that $\tau\le
M$, then
\begin{eqnarray*}
E\bigl\vert U_\varepsilon(\tau+\delta)-U_{\varepsilon}(\tau)
\bigr\vert^p&=& EE \bigl(\bigl\vert U_\varepsilon(\tau+\delta)-U_{\varepsilon}(
\tau )\bigr\vert^p|\F _{\tau}^{U_\eps} \bigr)
\\
&\le& C (M+\delta)^p \bigl(\delta\vee\delta^{{p}/ 2} \bigr),
\end{eqnarray*}
whenever
$1<p<1+\beta$ and $\varepsilon\le{a^{1+\beta}} \wedge\frac{t_0} {
M+\delta}$.\vspace*{1pt} This and the Markov inequality show that condition (ii)
[or equivalently, \eqref{e:Aldous2}] is also satisfied.

As already indicated, using Aldous' criterion we obtain the tightness
of the family $(U_{\varepsilon_n})_{n\geq1}$, which together with the
already proved convergence of finite dimensional distributions implies
that $(U_{\varepsilon_n})_n$ converges in law
to $(-t Z(t), t\ge0)$ with respect to the Skorokhod topology on
$D([0,\infty))$.

\textit{Step} 2.
For $b>0$, define
%
\begin{equation}
Z_{\varepsilon}^{(b)}(t)= \biggl(\frac{1} b
\ind_{[0,b]}(t)+\frac{1}t \ind _{(b,\infty)}(t)
\biggr)U_\varepsilon(t). \label{e:4.60}
\end{equation}
Recall that if $f\dvtx\R_+\mapsto\R$ is continuous, then the mapping
$w\mapsto fw$ is continuous from $D([0,\infty))$ into itself. Hence,
the result of step 1 implies that for any $b>0$, as $\varepsilon\to0$,
the family of processes $(Z^{(b)}_{\varepsilon})_{\varepsilon>0}$
converges in law to the process $Z^{(b)}$ defined by
\[
Z^{(b)}(t)= \frac{t} b\ind_{[0,b]}(t)Z(t)+
\ind_{(b,\infty)}(t)Z(t), \qquad t\geq0,
\]
with respect to the Skorokhod topology on $D([0,\infty))$.

\textit{Step} 3. We will next estimate the supremum norms of the
difference between $\tilde Y_\varepsilon$ and $Z^{(b)}_\varepsilon$,
and the difference between $Z$ and $Z^{(b)}$, respectively.
Fix any $1<p<1+\beta$ and suppose that $b\le t_0\wedge1$ and
$\varepsilon\le a^{1+\beta}$, where $t_0$ is as in \eqref{e:4.45} and
$a$ as in \eqref{e:A2}. Denote $\Vert f\Vert_\infty=\sup_{t\in\R
_+}\vert f(t)\vert$.

Using \eqref{e:tildeYe}--\eqref{e:U} and \eqref{e:4.60},
we have that
$\tilde Y_\varepsilon(t)-Z^{(b)}_\varepsilon(t) = U_\varepsilon
(t)(\frac{1}t - \frac{1}b)
\ind_{[0,b]}(t)$.
Therefore,
\[
\bigl\Vert\tilde Y_\varepsilon-Z^{(b)}_\varepsilon
\bigr\Vert_\infty\le\sup_{0\le
t\le b}\bigl\vert\tilde
Y_\varepsilon(t)\bigr\vert \le2\sup_{0\le t\le b}\bigl\vert
Y_\varepsilon(t)\bigr\vert,
\]
where \eqref{e:4.45} was used in the final estimate.
Lemmas \ref{lem:eqY} and \ref{L:bbl1} imply
\[
\sup_{0\le t\le b}\bigl\vert Y_\varepsilon(t)\bigr\vert\le2 \varepsilon
^{-{1}/{(1+\beta)}} \sup_{0\le t\le b}\bigl\vert M(\varepsilon t)\bigr\vert.
\]
Hence,
decomposing $M$ similarly as it was done for $U_\varepsilon$ in step 1 and
applying Doob's inequality for $M$, we obtain
%
\begin{eqnarray}
\label{e:4.61} &&E \bigl\Vert\tilde Y_\varepsilon-Z^{(b)}_\varepsilon
\bigr\Vert_\infty^p
\nonumber
\\[-8pt]
\\[-8pt]
\nonumber
&&\qquad \le C_1(p) \bigl(E\bigl\vert
\varepsilon^{-{1}/{(1+\beta)
}}M^{(1)}(\varepsilon b)\bigr\vert^p +E
\bigl\vert\varepsilon^{-{1}/{(1+\beta)}}M^{(2)}(\varepsilon b)\bigr\vert^p
\bigr),
\end{eqnarray}
where
\begin{eqnarray*}
M^{(1)}(\varepsilon b)&=&\int_0^{\varepsilon b}
\int_{[0,\varepsilon b]\times[0,\varepsilon^{{1}/{(1+\beta)}}]} y\hat\pi(du\,dy),
\\
M^{(2)}(\varepsilon b)&=& \int_{[0,\varepsilon b]\times
(\varepsilon^{{1}/{(1+\beta)}},1]} y\hat\pi(du\,dy).
\end{eqnarray*}
By mimicking the arguments of step 1, we obtain
\[
E\bigl\vert\varepsilon^{-{1}/{(1+\beta)}}M^{(1)}_{\varepsilon b}
\bigr\vert^2 \le C \varepsilon^{-{2}/{(1+\beta)}}\int_0^{\varepsilon b}
\int_0^{\varepsilon^{{1}/ {(1+\beta)}}} y^{-\beta}\,dy
\,du=C_1(p) b,
\]
and, relying on $\eps\ll\eps\eps^{{(p-1-\be)}/{(1+\be)}}=\eps
^{{p}/{(1+\be)}}$, we also obtain
\begin{eqnarray*}
&& E\bigl\vert\varepsilon^{-{1}/{(1+\beta)}}M^{(2)}_{\varepsilon b}
\bigr\vert^p \\
&&\qquad\le C_2(p) \varepsilon^{-{p}/{(1+\beta)}}\int
_0^{b \varepsilon} \biggl(\int_{\varepsilon^{{1}/{(1+\beta)}}}^a
y^{p-2-\beta}\,dy + \int_a^1
y^{p-2} \Lambda(dy) \biggr)\,du
\\
&&\qquad\le C_3(p) b.
\end{eqnarray*}
Together with \eqref{e:4.61} and Jensen's inequality, for $0<p<1+\beta
$, $b\le t_0\wedge1$ and $\varepsilon\le a^{1+\beta}$, this implies
%
\begin{equation}
\label{e:4.62} E\bigl \Vert\tilde Y_\varepsilon-Z^{(b)}_\varepsilon
\bigr\Vert_\infty^p\le C(p) b^{{p}/2},
\end{equation}
where $C(p)$ is some finite constant, uniform in $\varepsilon$.

For the processes $Z^{(b)}$ and $Z$, we again have
\[
E\bigl\Vert Z - Z^{(b)}\bigr\Vert_\infty\le\sup_{t\le b}
\bigl\vert Z(t)\bigr\vert.
\]
Since $Z$ is a solution of \eqref{e:eqZ}, we can again apply Lemma~\ref
{L:bbl1} and Doob's inequality to $L$,
a $(1+\beta)$-stable L\'evy process, to derive
%
\begin{equation}
\label{e:4.63} E\bigl\Vert Z - Z^{(b)}\bigr\Vert_\infty^p
\le C_1(p)E\bigl\vert L(b)\bigr\vert^p\le C_2(p)b^{{p}/{(1+\beta)}}
\end{equation}
for some $C_2(p)<\iy$.

\textit{Step} 4. Finally, we prove the convergence $\tilde
Y_\varepsilon
\Longrightarrow-Z$ as $\varepsilon\to0$.
Let $d^0_\infty$ denote the Skorokhod metric on $D([0,\infty))$ as
defined in \cite{Billingsley}, page 168. It is clear that $d^0_\infty
(f,g)\le\Vert f-g\Vert_\infty$ for any two $f,g\in D([0,\infty))$.

It suffices to show that, whenever $F\dvtx D([0,\infty))\mapsto
D([0,\infty
))$ is a given bounded and uniformly continuous function, we have
%
\begin{equation}
\lim_{\varepsilon\to0}\bigl\vert EF(\tilde Y_\varepsilon)-EF(Z)\bigr\vert=0.
\label{e:4.64}
\end{equation}
By the conclusion of step 2, for any $b>0$, we have $E|F(Z_\varepsilon
^{(b)})- EF(Z^{(b)})|\to0$. Hence, \eqref{e:4.64} follows by the
triangle inequality, the uniform continuity of $F$, estimates \eqref
{e:4.62} and \eqref{e:4.63} and the Markov inequality and the above
discussion. The argument based on addition and subtraction of
intermediate terms is standard, and the details are left to the reader.

\section{On robustness with respect to the choice of speed}
\label{S:robust}
\subsection{Proof of Theorem \texorpdfstring{\protect\ref{thmm:main a}}{1.4}}
Recall $\Psi$, $\Psi^*$ and $v$ defined in \eqref{e:Psi}, \eqref
{e:Psi_star} and~\eqref{e:vt}, respectively. Furthermore, recall that $v^*$\vadjust{\goodbreak} is defined
in terms of $\Psi^*$ as $v$ is defined in terms of $\Psi$.
Due to \eqref{e:5.4a},
one can easily see that
\[
\sup_{t\in[0,T]} \frac{1}{\eps^{{1}/{(1+\be)}}} \biggl\llvert
\frac{v_{\eps t}} {v_{\eps t^*} } -1 \biggr\rrvert =
O \bigl(\eps^{1-1/(\be+1)} \bigr)\qquad
\mbox{as $\eps\to0$}.
\]
Since
%
\begin{equation}\qquad\quad
\frac{1}{\eps^{{1}/{(1+\be)}}} \biggl( \frac{N_{\eps t}} {v_{\eps t}^* } -1 \biggr) =
\frac{1}{\eps^{{1}/{(1+\be)}}}
\biggl( \frac{N_{\eps t}} {v_{\eps t}} -1 \biggr) \times\frac{v_{\eps t}} {v_{\eps t} ^*} +
\frac{1}{\eps^{{1}/{(1+\be)}}} \biggl( \frac{v_{\eps t}} {v_{\eps t}^* } -1 \biggr), \label{e:8.1}
\end{equation}
one can conclude Theorem~\ref{thmm:main a}(a) directly from Theorem~\ref{thmm:main} and \eqref{e:2.4}.

We now turn to the proof of part (b).
Let us denote $w_t=K_1t^{-{1}/\beta}$ for $K_1$ from~\eqref{e:K1}.
Observe that an analogue of \eqref{e:8.1}, with $v^*$ replaced by $w$,
implies that it suffices to show
%
\begin{equation}
\label{e:8.2} \lim_{t\to0} t^{-{1}/{(1+\beta)} } \biggl(
\frac{v_t}{w_t}-1 \biggr)=0.
\end{equation}
Also note that $w$ is related to
$\Psi^{(\beta)}(q)=\frac{A\Gamma(1-\beta)}{\beta(1+\beta
)}q^{1+\beta}$
via relation
\[
t=\int_{w_t}\frac{1}{\Psi^{(\beta)}(q)}\,dq,
\]
the same way that $v$ is related to $\Psi$ [see \eqref{e:vt}]. Recall
that from \eqref{e:3.9} we already know $\lim_{q\to\infty}\Psi
(q)/\Psi^{(\beta)}
(q)=1$. We will need a more precise comparison of $\Psi$ and $\Psi
^{(\beta)}$.

Let $a\le\frac{1}2$ be such that $\Lambda$ has a density $g$ on $[0,a]$
satisfying \eqref{e:A2} and, moreover, $\vert y^\beta g(y)-A\vert\le C
y^\alpha$ on $[0,a]$. Such $a$ exists by the assumptions.

Observe that [similarly to derivation of \eqref{e:3.9}]
%
\begin{equation}
\Psi^{(\beta)}(q)=A q^2 \int_0^1
\int_0^r \int_0^\infty
e^{-qyu}y^{-\beta}\,dy \,du\,dr. \label{e:8.2a}
\end{equation}
Therefore, by Lemma~\ref{lem:3.1}(ii) and (iii),
we have
%
\begin{eqnarray}\label{e:8.3}
\qquad\Psi(q)=\Psi^*_a(q)+ O(q) =\Psi^{(\beta)}(q) +
R_1(q)- R_2(q)+O(q),
\nonumber
\\[-8pt]
\\[-8pt]
\eqntext{q\ge1,}
\end{eqnarray}
where
%
\begin{equation}
R_1(q)=q^2 \int_0^1
\int_0^r\int_0^a
e^{-qyu} \bigl(g(y)-Ay^{-\beta
} \bigr)\,dy\,du\,dr \label{e:R1}
\end{equation}
and
%
\begin{equation}
\label{e:R2} R_2(q)=q^2 A\int_0^1
\int_0^r \int_a^\infty
y^{-\beta}e^{-qyu}\,dy\,du\,dr.
\end{equation}
Due to the assumptions, we have
%
\begin{equation}
\label{e:R1 ad} \bigl|R_1(q)\bigr|\le C q^2\int
_0^1\int_0^r
\int_0^a y^{\alpha-\beta}e^{-qyu}\,dy
\,du\,dr.
\end{equation}
If $\alpha<\beta$, then (this is simpler than the proof of Lemma~\ref
{lem:3.2})
\[
\bigl|R_1(q)\bigr|\le C \Gamma(1+\alpha-\beta)q^{1+\beta-\alpha}\int
_0^1\int_0^r
u^{\beta-\alpha-1}\,du \,dr = O \bigl(q^{1+\beta-\alpha} \bigr).
\]
If $\alpha\ge\beta$, then by \eqref{e:R1 ad} we have
\begin{eqnarray*}
\bigl|R_1(q)\bigr| &\le& C q^2 a^{\alpha-\beta}\int
_0^1 \int_0^a
e^{-qyu}\,dy\,du
\\
&\le& Ca^{\alpha-\beta} \biggl(q^2\int_0^{{1}/q}
a \,du+q\int_{{1}/q}^1 \frac{1-e^{-qau}}{u}\,du
\biggr)
\\
&\le& Ca^{\alpha-\beta} \biggl(aq +q\int_{{1}/q}^1
\frac
{1}{u}\,du \biggr) =O \bigl(q(\log q+1) \bigr).
\end{eqnarray*}
For $R_2$, we have
%
\begin{equation}
R_2(q)\le Aq^2 \int_a^\infty
y^{-\beta} \int_0^1 e^{-qyu}\,du
\,dy \le Aq\int_a^\infty y^{-\beta-1}
\,dy=O(q). \label{e:8.5}
\end{equation}
Hence, from \eqref{e:8.3}, it follows that
%
\begin{equation}
\Psi(q)=\Psi^{(\beta)}(q) + O \bigl(q^{1+\beta-\alpha} \bigr) + O \bigl(q(
\log q+1) \bigr). \label{e:8.4}
\end{equation}
To prove \eqref{e:8.2},
we adapt the technique of Lemma~\ref{lem:2.2}(iii).
In particular, let us consider $v^{(n)}$ and $w^{(n)}$ defined by
\[
t=\int_{v^{(n)}_t}^n \frac{1}{\Psi(q)}\,dq \quad\mbox{and}\quad
t=\int_{w^{(n)}_t}^n \frac{1}{\Psi^{(\beta)}(q)}\,dq,
\]
and the following analogue of \eqref{E:last id}:
%
\begin{eqnarray}\label{e:8.6}
&&\log\frac{w_t^{(n)}}{v_t^{(n)}}+ \int_0^t \biggl[
\frac{\Psi
^{(\beta)}
(w_s^{(n)})}{w_s^{(n)}} - \frac{\Psi^{(\beta)}
(v_s^{(n)})}{v_s^{(n)}} \biggr] \,ds
\nonumber
\\[-8pt]
\\[-8pt]
\nonumber
&&\qquad= \int_0^t
\frac{\Psi(v_s^{(n)})-\Psi^{(\beta)}
(v_s^{(n)})}{v_s^{(n)}} \,ds
\end{eqnarray}
(note that if $n\ge2$ and $t$ is sufficiently small, then
$w_s^{(n)}\ge1$ for $s\le t$). Also, observe that
$v_s^{(n)}\nearrow v_s$, $w_s^{(n)}\nearrow w_s$ as $n\to\infty$.
Lemma~\ref{L:bbl1} implies that for sufficiently small $t\le t_0$
(with $t_0$ uniform in $n\ge2$) we have
\[
\biggl\vert\log\frac{w_t^{(n)}}{v_t^{(n)}}\biggr\vert\le2 \int_0^t
\frac
{\vert\Psi (v_s^{(n)})-\Psi^{(\beta)}(v_s^{(n)})\vert
}{v_s^{(n)}}\,ds.
\]
Using \eqref{e:8.4} and $v_s^{(n)}\le v_s\le C s^{-{1}/\beta}$ for
small $s$ [see \eqref{e:3.10a}], we obtain
%
\begin{eqnarray} \label{e:8.7}
\biggl\vert\log\frac{w_t^{(n)}}{v_t^{(n)}}\biggr\vert&\le& C \biggl( \int_0^t
(v_s)^{\beta
-\alpha}\,ds +\int_0^t
\log(v_s)\,ds \biggr)
\nonumber
\\[-8pt]
\\[-8pt]
\nonumber
&=&O \bigl(t^{\alpha/\beta} \bigr) +O \biggl(t
\log \frac{1}t \biggr).
\end{eqnarray}
Letting $n\to\infty$, we see that the same estimate holds also for $
\vert\log\frac{w_t}{v_t}\vert=\vert\log\frac{v_t}{w_t}\vert$.
In particular, $\lim_{t\to0+}\log\frac{v_t}{w_t}=0$, and so
$\vert\frac{v_t}{w_t}-1\vert\sim\vert\log\frac{w_t}{v_t} \vert$
for small $t$.
We conclude that \eqref{e:8.2} holds since $\frac{\alpha}{\beta
}>\frac
{1}{1+\beta}$, completing the proof.

\subsection{Limitations of robustness}
\label{S:robust b}
In this section, we provide an instructive counterexample,
announced in both the \hyperref[sec1]{Introduction} and Remark~\ref
{rem:natural}.
A careful reader will note that the just made arguments proving Theorem~\ref{thmm:main a} are close to optimal, in that the power $\alpha
=\frac
{\be}{1+\be}$ should be critical for \eqref{e:8.2}.
Without making any general statements to this end, let us fix $\alpha
\in(0,\frac{\be}{1+\be})$ and consider $\La$ such that
\[
\La(dy) = g(y)\,dy,\qquad y \in[0,1], \mbox{ where } g(y):=y^{-\be} \bigl(1 +
y^\alpha \bigr), y\in(0,1].
\]
We keep the notation of the previous section, setting $A=1$ (note that
hence $\La$ is not anymore a probability measure but, as mentioned in
the second paragraph of the \hyperref[sec1]{Introduction}, all our
results continue to
hold with appropriately modified constants). In particular, $v$ and $w$
are as in \eqref{e:8.2}, up to the same positive multiple.
We will show that
%
\begin{equation}
\label{E:show two} {t^{-{1}/{(1+\be)}}} \biggl(\frac{w_{t}}{v_{t}} -1
\biggr) \qquad
\mbox{is unbounded as } t\to0,
\end{equation}
and that therefore the statement of Theorem~\ref{thmm:main a}(b)
cannot hold in this particular case.

As in \eqref{e:8.3} and \eqref{e:8.5} (with $a=\frac{1}2$), we have
\[
\Psi(q) - \Psi^{(\beta)}(q) = R_1(q) +O(q), \qquad q\ge1.
\]
Now $R_1$ can be written explicitly as
\[
R_1(q)= q^2 \int_0^1
\int_0^r\int_0^{{1}/2}
e^{-qyu}y^{\alpha
-\beta}\,dy\,du\,dr.
\]
Note that $R_1$ is again of the form \eqref{e:Psi_star} where $\La$ is
given by $\Lambda_{\beta-\alpha}(dy)=y^{\alpha-\beta}\ind
_{[0,{1}/2]}(y)\,dy$.
By \eqref{e:8.2a}, \eqref{e:R2} and \eqref{e:8.5} with $\beta$ replaced
by $\beta-\alpha$, we obtain
%
\begin{eqnarray} \label{e:Psib_2}
{\Psi(q)-\Psi^{(\beta)}(q)}&=&\Psi^{(\beta-\alpha)}(q)+O(q)
\nonumber
\\[-8pt]
\\[-8pt]
\nonumber
&=& D
q^{1+\beta
-\alpha} +O(q),\qquad q\ge1,
\end{eqnarray}
where $D$ is a positive constant that can be written explicitly.

Recall the expression for $\Psi^{(\beta)}$ given just after \eqref{e:8.2}.
It is easy to check that one can let $n\to\infty$ in \eqref{e:8.6},
and obtain
%
\begin{equation}
\label{e:8.8} \log\frac{w_t}{v_t}+C\int_0^t
\bigl(w_s^\beta- v_s^\beta \bigr)
\,ds=D \int_0^t v_s^{\beta-\alpha}
\,ds +O(t)
\end{equation}
for all sufficiently small $t$, where $C$ and $D$ are positive
constants (their exact value is not important for our purposes).
As usual, this is done via uniform (in small $t$ and in $n$) control of
the RHS in \eqref{e:8.6}; see \eqref{e:8.7} for a similar argument.
By
\eqref{e:3.10a}, it follows that
%
\begin{equation}
\label{e:8.10} \int_0^t v_s^{\beta-\alpha}
\,ds \sim C_1 t^{{\alpha}/ {\beta}}.
\end{equation}

Let us suppose that the function given in \eqref{E:show two} is bounded
near $0$.
Since $\alpha<\frac\beta{1+\beta}$, this implies that
\[
\biggl\vert\frac{w_t}{v_t}-1\biggr\vert=o \bigl(t^{\alpha/\beta} \bigr)\qquad
\mbox{as }t
\to0,
\]
hence also
%
\begin{equation}
\biggl\vert\log\frac{w_t}{v_t}\biggr\vert \vee\biggl\vert\frac{v_t}{w_t}-1\biggr\vert=o
\bigl(t^{\alpha/\beta} \bigr)\qquad \mbox{as }t\to0. \label{e:assum o}
\end{equation}
By an elementary application of Taylor's formula, we have
\[
\bigl\vert w_s^\beta-v_s^\beta\bigr\vert=
\biggl\vert1- \biggl(\frac
{v_s}{w_s} \biggr)^\beta \biggr\vert
w_s^\beta\sim\beta\biggl\vert1-\frac
{v_s}{w_s}\biggr\vert
w_s^\beta\qquad \mbox{as } s\to0,
\]
and since $w_s^\be= K_1^\be s^{-1}$, we conclude
\begin{eqnarray*}
\int_0^t \bigl\vert{w_s}^\be-
{v_s}^\be\bigr\vert \,ds &\le& C \beta K_1^\be
\int_0^t \frac{1}{s}\biggl\vert1-
\frac{v_s}{w_s}\biggr\vert \,ds
\\
&=&C \beta K_1^\be\int_0^t
o \bigl(s^{-1+{\alpha}/{\be}} \bigr) \,ds = o \bigl(t^{ {\alpha} /{\be}}
\bigr).
\end{eqnarray*}
This together with \eqref{e:assum o} is in clear contradiction with
\eqref{e:8.10} and \eqref{e:8.8}.
We conclude that the opposite of \eqref{e:8.2} must hold, or
equivalently, that there must exist a positive constant $c$ and a
sequence of times $(t_n)_n$ such that $t_n\to0$ and
\[
\biggl\llvert \frac{v_{t_n}}{w_{t_n}} -1 \biggr\rrvert \geq c
 (t_n)^{{\alpha}/{\be}},
\]
and joint with $\alpha\in(0,\frac{\be}{1+\be})$, this easily implies
\eqref{E:show two}.

\section*{Acknowledgment} We would like to thank the anonymous
referee for a careful reading of the paper, and for several helpful
suggestions that improved the presentation.

%





\printaddresses

\end{document}